%% file: eigenvalue-arxiv.tex
\documentclass{jssc}

\def\textsubscript#1%
{$_{\text{#1}}$}
\newcommand*\supercite[1]{\textsuperscript{\cite{#1}}}
\def\cdd{\mbox{\boldmath$\cdot$}~}

\input{amssym.def}
\include{graphix}

\usepackage{mathtools}
\usepackage{amsmath}
\usepackage{amsfonts}
\usepackage{amssymb}
\usepackage{mathrsfs}

\usepackage{lineno}

\usepackage{indentfirst} 
\usepackage{float} 
\usepackage{color} 
\usepackage{graphicx} 

\usepackage{enumerate} 

\usepackage{algorithm} 

\usepackage{multirow} 

\usepackage{svg} 

\usepackage{datetime} 

\newcommand{\mathbbR}{\mathbb{R}}

\newcommand{\rme}{${\rm e}$}
\newcommand{\itOmega}{{\it \Omega}}

\input jsscN.tex
\begin{document}

\title{A Deep Learning Method for Computing Eigenvalues of the Fractional Schrödinger Operator$^*$}
{\uppercase{Guo} Yixiao \cdd 
	\uppercase{Ming} Pingbing}
{
\uppercase{Guo} Yixiao \cdd \uppercase{Ming} Pingbing (Corresponding author) \\
{\it LSEC, Institute of Computational Mathematics and Scientific/Engineering Computing, AMSS, Chinese Academy of Sciences, No. 55, East Road Zhong-Guan-Cun, Beijing 100190, China; School of Mathematical Sciences, University of Chinese Academy of Sciences, Beijing 100049, China.}
Email: guoyixiao@lsec.cc.ac.cn; mpb@lsec.cc.ac.cn.\\
} 
{$^*$ This work is supported by National Natural Science Foundation of China through Grant No. 11971467.\\
{$^\diamond$}}

\drd{DOI: }{Received: x x 20xx}{ / Revised: x x 20xx}


\dshm{20XX}{XX}{\uppercase{deep learning for eigenvalues of fractional Schrödinger}}{
	\uppercase{Guo Yixiao} $\cdd$ 
	\uppercase{Ming Pingbing}}

\Abstract{
We present a novel deep learning method for computing eigenvalues of the fractional Schrödinger operator. 
Our approach combines a newly developed loss function with an innovative neural network architecture that incorporates prior knowledge of the problem. 
These improvements enable our method to handle both high-dimensional problems and problems posed on irregular bounded domains.
We successfully compute up to the first $30$ eigenvalues for various fractional Schrödinger operators. 
As an application, we share a conjecture to the fractional order isospectral problem that has not yet been studied.}      




\section{Introduction}

Fractional partial differential equations have proven effective in modeling anomalous diffusion phenomena\supercite{intro_ad1_metzler2000random}, turbulent flows\supercite{intro_tf1_shlesinger1987levy}, porous media flows\supercite{intro_pmf1_de2011fractional} and many others. 
In the field of quantum mechanics, Laskin introduced the fractional Schrödinger equation\supercite{intro_qm1_laskin2000fractional, intro_qm2_laskin2002fractional} by using Levy flights to replace the classical Brownian motion. 
This equation has been used to reveal some novel phenomena\supercite{intro_qm3_longhi2015fractional, intro_qm4_zhang2016pt}, which could not be explained by the standard Schrödinger equation.

Because of the non-locality and the singularity, most of the numerical methods for solving  fractional partial differential equations focus on one-dimensional problems. For two-dimensional problems, the most popular method is an adaptive finite element method proposed by Ainsworth and Glusa\supercite{intro_AFEM1_ainsworth2017aspects, intro_AFEM2_ainsworth2018towards}, while \cite{intro_FD2_del2018robust}
also presents some modified finite difference methods. 
In addition, the authors in \cite{intro_Spectral1_mao2016efficient, intro_Spectral2_xu2018spectral} studied the spectral methods to solve the fractional partial differential equations, but applying them in irregular domains presents significant challenges.
Another innovative method is the walk-on-spheres (WOS) method\supercite{intro_WOS1_kyprianou2018unbiased, intro_WOS2_shardlow2019walk} and its extension\supercite{intro_WOS3_sheng2022efficient}.
Those methods use the Feynman-Kac formula to calculate the solution at any given point by simulating a large number of paths of the $2s$-stable Lévy process.
We refer to \cite{du2020acta, Karniadakis2020whatis, intro_num_method3_bonito2018numerical} for a comprehensive review of the numerical methods for fractional partial differential equations.

In this paper, we study the eigenvalue problem of the fractional Schrödinger operator in a bounded domain.
\begin{equation}
\left \{
\begin{aligned}
& (-\Delta)^s u(x) + V(x) u(x) = \lambda u(x), & \ x \in \itOmega \subset \mathbbR^d, \\
& u(x) = 0, & \ x \in \itOmega^c \triangleq \mathbbR^d \backslash \itOmega.
\end{aligned}\right.
\label{problem_schrodinger}
\end{equation}
Here, $V(x)$ is a real-valued potential energy function. The precise definition of the operator $(-\Delta)^s$ with $s\in(0,1)$ will be presented in the subsequent section. When $V(x) = 0$, the problem simplifies to the eigenvalue problem of the fractional Laplace operator:
\begin{equation}
\left \{
\begin{aligned}
& (-\Delta)^s u(x) = \lambda u(x), & \ x \in \itOmega, \\
& u(x) = 0, & \ x \in \itOmega^c.
\end{aligned}\right.
\label{problem_laplace}
\end{equation}
To the best of our knowledge, the exact eigenvalue for any of the aforementioned problems is still unknown, and no analytical formulas even exist for calculating them under any specific  circumstance. However, several numerical methods have been developed to tackle this challenge.

In \cite{dyda2017eigenvalues, dyda2012fractional}, the author presents a spectral method capable of computing very tight bounds for the eigenvalues of the fractional Laplacian within a unit ball in any dimension.
Additionally, \cite{bao2020jacobi} introduces a spectral method that approximates the eigenvalues of the one-dimensional fractional Schrödinger operator. This method employs the same basis functions as \cite{intro_Spectral1_mao2016efficient} and can be extended to handle problems in hypercubes. 
Furthermore, a finite element method is displayed in \cite{FEM_borthagaray2018finite} for approximating the eigenvalues within arbitrary bounded domains in one and two dimensions.

The use of neural networks for computing eigenvalues and eigenvectors dates back to the 1990s\supercite{intro_NN_early1_samardzija1991neural, intro_NN_early2_cichocki1992neural}. 
While the recent efforts focus on solving the eigenvalue problem of many-body quantum systems\supercite{intro_NN_mb1_carleo2017solving, intro_NN_mb2_choo2020fermionic, intro_NN_mb5_han2019solving}. These methods use neural networks to represent the underlying functions and incorporate techniques such as the variational Monte Carlo method and stochastic optimizers to approach the ground state of the quantum system.
\cite{intro_NN_eigen1_han2020solving, intro_NN_eigen2_li2022semigroup, intro_NN_eigen3_simonnet2022deep, intro_NN_eigen5_zhang2022solving, intro_NN_eigen9_elhamod2022cophy, intro_NN_eigen10_finol2019deep} also exploit neural networks to solve eigenvalue problems. 
However, there is little literature on the numerical analysis of solving eigenvalue problems with neural networks except\supercite{intro_NN_eigen_theory_lu2022priori}.

In this work, we propose a deep learning method to solve the eigenvalue problem associated with the fractional Schrödinger operator. 
We convert the original eigenvalue problem into a sequence of minimization problems with constraints. By solving these minimization problems sequentially, we can determine the eigenmodes in ascending order of their eigenvalues.
In addition, we introduce a novel loss function for solving these minimization problems. This loss function incorporates penalty terms to efficiently handle the orthogonal constraints.
Moreover, we design a new neural network architecture that incorporates various feature functions. These feature functions are derived from prior knowledge of the underlying problem, particularly focusing on aspects such as singularity and boundary conditions.

To evaluate the accuracy of our method, we firstly solve a range of examples that could be addressed by the spectral methods also. The dimension of these problems varies from $1$ to $9$ and our network employs approximately 3,500 parameters. 
The relative error of our method is less than $0.1\%$ for the first $5$ eigenvalues. We further calculate up to $30$ eigenvalues while maintaining an error of less than $1\%$. 
Subsequently, we test our method in general domains where the spectral methods could not be used and compare the results with those obtained by the finite element method. The results demonstrate that our method produces more accurate results than the finite element method over the finest mesh. 

We also implement our method for a fractional Schrödinger operator with an inverse square potential function in three dimensions. 
By computing the first $30$ eigenvalues, we observe that the order of the eigenvalues exchanges as the fractional order varies in this example.
Additionally, we provide our estimations of eigenvalues for problems with different potential functions that have never been tested and exhibit eigenfunctions for various problems. All these examples collectively demonstrate the accuracy and efficiency of our method.

As an application, we apply the method to the fractional version of the isospectral problem. Based on our numerical result, we conjecture that even if the spectra of two domains are identical for the Laplacian, they would not be the same for the fractional Laplacian.

The rest of the paper is as follows. In Section \ref{Section_Formulation}, we begin with a discussion about the fractional operator and present the newly devised loss function. In Section \ref{Section_Deep_learning_scheme}, we display a deep learning scheme for solving the eigenvalue problem with a novel neural network architecture. We show the results of numerous numerical experiments in Section \ref{Section_Numerical_results} and compare them with other existing methods. We apply our method to the fractional version of the isospectral problem in Section \ref{Section_Isospectral} and draw the conclusions in Section \ref{Section_Conclusion}.

\section{Formulation}\label{Section_Formulation}
The fractional Laplace operator $(-\Delta)^s$ has many different definitions\supercite{Karniadakis2020whatis, kwasnicki2017ten} and these definitions are not equivalent in all circumstances. In this paper, we adopt the Ritz definition as follows.
\begin{equation}
(-\Delta)^s u(x) = C_{d,s} \int_{\mathbbR^d} \frac{u(x)-u(y)}{\|y-x\|^{d+2s}} dy \quad \text{for all} \ x \in \mathbbR^d,
\label{integral_laplacian}
\end{equation}
where
\begin{equation}
C_{d,s}{:}= \frac{2^{2s}s\Gamma(s+d/2)}{\pi^{d/2}\Gamma(1-s)}
\end{equation}
and $\|y-x\|$ represents the distance between $x$ and $y$, i.e., $\|y-x\|=\|y-x\|_{\ell_2}$. The integral in \eqref{integral_laplacian} should be interpreted in the principal value sense. This definition is also known as the integral definition in the literature and it is equivalent to the definition via the Fourier transform\supercite{valdinoci2009long}
\begin{equation}
(-\Delta)^s u(x) = \mathcal{F}^{-1}\{|\xi|^{2s}\mathcal{F}\{u\}(\xi)\}(x),
\label{fourier_laplacian}
\end{equation}
where $\mathcal{F}$ and $\mathcal{F}^{-1}$ represent the Fourier transform and the inverse Fourier transform, respectively.

To discuss the variational form of the eigenvalue problem \eqref{problem_schrodinger}, we first define the standard fractional order Sobolev space as 
\begin{equation}
H^s(\mathbbR^d){:}= \left\{u \in L^2(\mathbbR^d) : \int_{\mathbbR^d} (1+|\xi|^{2s}) |\mathcal{F}u(\xi)|^2 d\xi < \infty\right\} .
\end{equation}
With the equivalent definition of the fractional Laplacian, it can be expressed as
\begin{equation}
H^s(\mathbbR^d){:}= \left\{u \in L^2(\mathbbR^d) : \|u\|_{H^s(\mathbbR^d)} < \infty\right\},
\end{equation}
where the norm is
\begin{equation}
\|u\|_{H^s(\mathbbR^d)}^2 = \|u\|^2_{L^2(\mathbbR^d)} + |u|^2_{H^s(\mathbbR^d)}
\end{equation}
with the seminorm is
\begin{equation}
|u|_{H^s(\mathbbR^d)}^2 = \int_{\mathbbR^d} \int_{\mathbbR^d} \frac{[u(y)-u(x)]^2}{\|x-y\|^{d+2s}} dy dx.
\end{equation}

For the eigenvalue problem with homogeneous Dirichlet boundary condition, we find a solution in 
\begin{equation}
H^s_c(\itOmega){:}= \left\{u \in H^s(\mathbbR^d) : u(x) = 0 \text{ for all } x \in \itOmega^c \right\}.
\end{equation}
A bilinear form associated with the fractional Schrödinger operator is derived from the nonlocal Green's first identity \cite[(1.22)]{du2020acta}
\begin{equation}
\begin{aligned}
a(u,v) {:}= & \ \int_{\itOmega} ((-\Delta)^s u(x) + V(x)u(x)) v(x) dx \\
= & \  \frac{C_{d,s}}{2} \int_{\mathbbR^d}\int_{\mathbbR^d} \frac{(u(x)-u(y))(v(x)-v(y))}{\|x-y\|^{d+2s}} dydx
+ \int_{\itOmega} V(x)u(x)v(x) dx.
\end{aligned}
\end{equation}

Using this bilinear form and the variational principle, the $k$th smallest eigenvalue is given by
\begin{equation}
\lambda_k = \min_{E} \max_{u \in E \backslash \{0\}}  
\frac{a(u,u)}{\|u\|_{L^2(\itOmega)}^2},
\end{equation}
where $E$ is a $k$-dimensional subspace of $H^s_c(\itOmega)$. 
\cite{property1_frank2016eigenvalue, property2_chen2005two, property3_bao2018fundamental} display some analytical results for the eigenvalue of the fractional Schrödinger operator.

Numerically solving the min-max problem to derive the eigenvalues is challenging and requires significant effort.
Therefore, we propose a new formulation that provides a more convenient approach to solve the problem by solely minimizing a loss function. Given the first $k$ eigenmodes $(\lambda_1, u_1), \cdots, (\lambda_k, u_k)$, the subsequent eigenvalue can be characterized as
\begin{equation}
\lambda_{k+1} = \min_{u \in E^{(k)}} \frac{a(u,u)}{\|u\|_{L^2(\itOmega)}^2},
\end{equation}
where
\begin{equation}
E^{(k)} = \{u \in H^s_c(\itOmega) \backslash \{0\}: u \perp u_i, 1 \leq i \leq k \}.
\end{equation}
We use a neural network to approximate the next eigenfunction and construct a loss function by incorporating penalty terms to handle the orthogonal constraints.
The loss function for computing the $k$th eigenvalue is
\begin{equation}\label{Loss}
L_{k}(u) = \frac{a(u,u)}{\|u\|_{L^2(\itOmega)}^2} + \beta \sum_{j=1}^{k-1}
\frac{(u,u_j)^2}{\|u\|_{L^2(\itOmega)}^2 \|u_j\|_{L^2(\itOmega)}^2},
\end{equation}
where $\beta$ is a penalty parameter and its value should be greater than $\lambda_k - \lambda_1$.
When $k=1$, the loss function reduces to
\begin{equation}
L_1(u) = \frac{a(u,u)}{\|u\|_{L^2(\itOmega)}^2},
\label{Loss1}
\end{equation}
and the penalty term becomes unnecessary.
We remark that the loss function \eqref{Loss} originates from the variational principle and does not rely on any specific properties of the fractional differential operators. Thus, this loss function can be employed for the standard Schrödinger operator and other ordinary differential operators that admit Dirichlet forms.

We use the deep learning scheme introduced in the next section to minimize the loss. After minimizing the loss, the neural network provides an approximation of the eigenfunction $\hat{u}_k$, while the eigenvalue $\lambda_k$ is approximated as
\begin{equation}
\hat{\lambda}_k = \frac{a(\hat{u}_k,\hat{u}_k)}{\|\hat{u}_k\|_{L^2(\itOmega)}^2}.
\label{MC_eigen}
\end{equation}
We iterate this process to obtain the eigenmodes one by one. Since the exact eigenfunction is unknown, the approximating function $\hat{u}_k$ will be used as a substitute in \eqref{Loss}.
Consequently, any numerical errors from previous calculations can impact the accuracy of eigenmodes with higher eigenvalues.
As the training progresses, the numerical error becomes larger and larger.
However, subsequent experiments demonstrate that our method can calculate the first dozens of eigenmodes with great precision.

\section{Deep Learning Scheme}\label{Section_Deep_learning_scheme}
In this section, we present our deep learning scheme for solving the eigenvalue problem \eqref{problem_schrodinger}. We introduce a new architecture and describe the Monte Carlo sampling method used to calculate the loss.

When applying deep learning to solve PDE problems, the fully connected neural network (FCNN) and the residual network (ResNet) are commonly employed to approximate functions. However, these architectures have some weaknesses.
They could not ensure the boundary conditions, and penalty terms are often added to the loss function to constrain the functions in general.
Moreover, these architectures tend to underperform when the solutions contain singular terms. 
According to \cite{boundary_order_grubb2015fractional}, solutions of fractional partial differential equations exhibit an $s$-order singularity near the boundary, i.e.,
\begin{equation}
u(x) \approx \text{dist}(x, \partial \itOmega)^s + v(x), \quad x \in \bar{\itOmega}
\end{equation}
where $v(x)$ is a smooth function over $\bar{\itOmega}$.
Furthermore, \cite{bao2020jacobi} also conjectures that the eigenfunctions of the fractional Schrödinger operator exhibit an $s$-order singularity near the boundary. 
To overcome these drawbacks, we design a new architecture that incorporates prior knowledge about the singular term of the functions and the boundary conditions.

Our architecture is based on FCNN, and a similar one is shown in \cite{architecture_khoo2019solving}. 
It consists of several hidden layers and an output layer, while we mainly change the formulation of the output layer. Each hidden layer is a fully connected layer with the same width. 
We denote the number of hidden layers as $l$ and the width of each layer as $m$.
The activation function we used is $\sigma = \tanh$ and we denote the vectorized function of $\sigma(x)$ as $\phi(x)$, i.e.,
\[
\phi(x) = (\sigma(x_1), \sigma(x_2), \dots, \sigma(x_m)).
\]

\begin{figure}
	\begin{center}		
		\includegraphics[width=0.6\columnwidth]{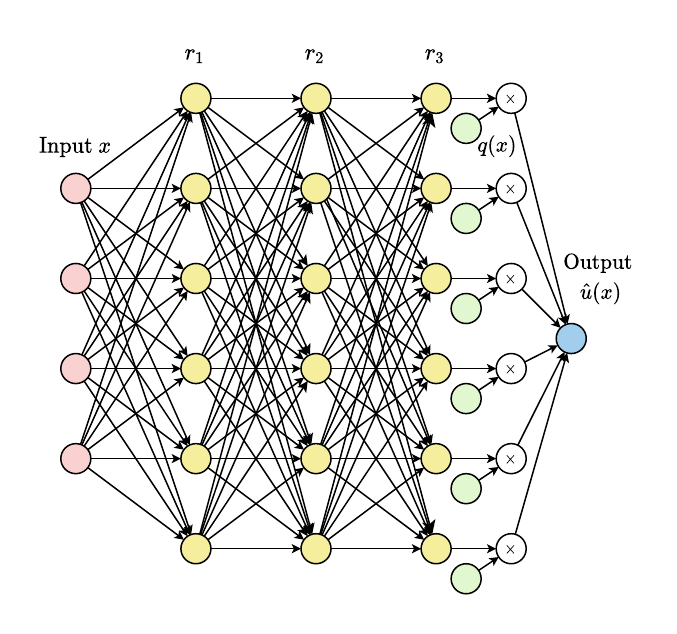}
		
		\small {\bf Figure 1}\ \ The neural network architecuture combined with feature funtcions.
	\end{center}
\end{figure}

The input layer contains a linear transformation that maps input values from $\mathbbR^d$ to $\mathbbR^m$, followed by an activation introducing non-linearity. The output of the first layer is given by
\begin{equation}
r_1 = \phi(W_1 x + b_1),
\end{equation}
where $W_1 \in \mathbbR^{m \times d}$ and $b_1 \in \mathbbR^{m}$. 
The other hidden layers also contain similar transformations, mapping values from $\mathbbR^m$ to $\mathbbR^m$, and the output of the $i$th layer is represented as
\begin{equation}
r_i = \phi(W_ir_{i-1} + b_i), \quad 2 \leq i \leq l,
\end{equation}
where $W_i \in \mathbbR^{m \times m}$ and $b_i \in \mathbbR^m$.
The output layer combines the outputs of the $l$th layer and a series of feature functions.
The final output of the entire network is given by
\begin{equation}
\hat{u}(x) = \sum_{j=1}^{m} W_{l+1}^{(j)} \cdot r_l^{(j)} \cdot q_j(x),
\end{equation}
where $W_{l+1} \in \mathbbR^{m}$, $W_{l+1}^{(j)}$ and $r_l^{(j)}$ is the $j$th component of $W_{l+1}$ and $r_l$, respectively. 
The feature functions are the key part of our method and play a critical role in ensuring the boundary conditions and capturing the singularity of the eigenfunctions.
The selection of these functions is based on prior knowledge and significantly influences the accuracy and efficiency of our method. We enforce that the feature functions belong to $H^s_c(\itOmega)$ to ensure the final output resides in the appropriate function space.
The choice of the feature functions for different examples will be provided in the next section.

The set of all parameters is defined as
\begin{equation}
\theta{:}= \{W_1, \dots, W_{l+1}, b_1, \dots, b_l\}.
\end{equation}
During each epoch, we calculate the loss function $L_k[u(x;\theta)]$ and employ stochastic optimization methods, such as the adaptive moment estimation (ADAM) optimizer, to update the parameters. 
After multiple epochs, the loss function will diminish to a value small enough and we will derive the approximated eigenmode.
To balance accuracy and efficiency, we gradually reduce the learning rate and increase the number of sampling points as the training progresses. 

Next, we present the details of the sampling technique used to estimate the loss function.
The calculation of $(u,u_j)$ and $\|u\|^2_{L^2(\itOmega)}$ is straightforward. We uniformly sample $N$ points $x_1, x_2, \dots, x_N$ in $\itOmega$ and calculate these values unbiasedly by
\begin{equation}
(u,u_j) \approx \frac{|\itOmega|}{N} \sum_{i=1}^N u(x_i) u_j(x_i), 
\end{equation}
and
\begin{equation}
\|u\|^2_{L^2(\itOmega)} \approx \frac{|\itOmega|}{N} \sum_{i=1}^N u(x_i)^2. 
\end{equation}

However, calculating the quadratic form $a(u,u)$ requires a more sophisticated approach since it is a double integral over $\mathbbR^d$ and the denominator becomes extremely small when the two variables are in close proximity. To address these difficulties, we pick a convex domain $D$ that contains $\itOmega$ and separate the integral into two parts.

\begin{equation}
\begin{aligned}
a(u,u){:}= \ & \frac{C_{d,s}}{2} \int_{\mathbbR^d} \int_{\mathbbR^d} \frac{[u(y)-u(x)]^2}{\|x-y\|^{d+2s}} dy dx \\
= \ & \frac{C_{d,s}}{2} \int_{D} \int_{D} \frac{[u(y)-u(x)]^2}{\|x-y\|^{d+2s}} dy dx + 
C_{d,s} \int_{D} \int_{D^c} \frac{[u(y)-u(x)]^2}{\|x-y\|^{d+2s}} dy dx \\
\triangleq \ & A_1 + A_2. \\
\end{aligned}
\end{equation}

For the first part $A_1$, the integral domain is bounded, but its value approaches infinity as $y$ approaches $x$. We alleviate this by applying a coordinate transformation to change the integral formulation
\begin{equation}
\begin{aligned}
A_1 = \ & \frac{C_{d,s}}{2} \int_{D} \int_{D} \frac{[u(y)-u(x)]^2}{\|x-y\|^{d+2s}} dy dx \\
= \ & \frac{C_{d,s}}{2} \int_{D} \int_{S^{d-1}} \int_0^{+\infty} \mathbf{1}_{x+w\xi \in D} \frac{[u(x+w\xi)-u(x)]^2}{w^{1+2s}} dwd\xi dx \\
= \ & \frac{C_{d,s}}{2} \int_{D} \int_{S^{d-1}} \int_0^{w^+}  \bigg[\frac{u(x+w\xi)-u(x)}{w}\bigg]^2 \frac{dw}{w^{2s-1}} d\xi dx.
\end{aligned}
\label{eqA1}
\end{equation}
Here, $S^{d-1}$ is the unit $(d-1)$-dimensional sphere and $w^+$ denotes the distance from the point $x$ to the boundary $\partial D$ along the direction $\xi$. The shape of $D$ should facilitate this calculation and its size should be as small as possible to enhance sampling efficiency. It is crucial that $D$ is a convex domain, which guarantees that any ray starting from $x$ intersects the boundary $\partial D$ only once. This makes the last equal sign in \eqref{eqA1} and simplifies the calculation.

For the numerical implementation, we first uniformly sample $x_1, x_2, \dots, x_N$ in $D$, and uniformly sample $\xi_1, \xi_2, \dots, \xi_N$ in $S^{d-1}$. Then, we calculate $w^+_i$ which depends on $x_i$ and $\xi_i$, and sample $w_i$ with the probability
\begin{equation}
P(w_i = w) = \mathbf{1}_{0 < w < w_i^+} \frac{1}{w^{2s-1}} \frac{2-2s}{(w_i^+)^{2-2s}}.
\end{equation}
To reduce numerical instability caused by dividing a very small amount, we set a minimum value $w_c$ for $w$. In our numerical experiments, we take $w_c = 10^{-4}$ and
\begin{equation}
\tilde{w}_i = \max(w_i, w_c) = \max(w_i, 10^{-4}).
\end{equation}
Through these steps, we derive a practicable and efficient sampling method for calculating $A_1$,
\begin{equation}
A_1 \approx \frac{C_{d,s}}{2N} \ |D| \ |S^{d-1}| \ \sum_{i=1}^N 
\bigg[ \big[\frac{u(x_i+\tilde{w}_i\xi_i)-u(x_i)}{\tilde{w}_i}\big]^2 \frac{(w^+_i)^{2-2s}}{2-2s}\bigg].
\end{equation} 

The calculation of $A_2$ is comparatively simpler, as $u(y) = 0$ when $y \in \itOmega^c$. We simplify the equation by calculating part of the integral in advance and this prevents sampling in an infinite domain.
\begin{equation}
\begin{aligned}
A_2 = \ & C_{d,s} \int_{D} \int_{D^c} \frac{[u(y)-u(x)]^2}{\|x-y\|^{d+2s}} dy dx \\
= \ & C_{d,s} \int_{D} \int_{S^{d-1}} \int_{w^+}^{+\infty} \frac{u(x)^2}{w^{1+2s}} dw d\xi dx \\
= \ & C_{d,s} \int_{D} \int_{S^{d-1}} u(x)^2 \frac{1}{(w^+)^{2s}}\frac{1}{2s} d\xi dx. \\
\end{aligned}
\end{equation}
By employing the same method to uniformly sample $x_i$ and $\xi_i$ in $D$ and $S^{d-1}$, we can approximate $A_2$ unbiasedly by
\begin{equation}
A_2 \approx \frac{C_{d,s}}{N} \ |D| \ |S^{d-1}| \ \sum_{i=1}^N 
\bigg[ u(x_i)^2 \frac{1}{(w^+_i)^{2s}}\frac{1}{2s}\bigg].
\end{equation}

\newpage
\section{Numerical results}
\label{Section_Numerical_results}
In this section, we present the numerical results for various problems. 
Firstly, we introduce the training configuration. We employ the networks introduced in the last section. Each of them contains $3$ hidden layers with a width of $40$, unless otherwise stated. 
The total number of parameters is approximately $3,500$. For each eigenvalue, we conduct $120,000$ epochs. 
Initially, the learning rate is set to $5\rme\text{-}03$, and $1,000$ points are sampled. After every $20,000$ epochs, we reduce the learning rate to one-fourth of its current value, and double the number of sampling points. 
To identify the $k$th eigenvalue $(k \geq 2)$, we set the penalty parameter $\beta$ to $4$ times the maximum eigenvalue we have found. This selection is sufficiently large to discover the subsequent correct eigenvalue. 
Using a too large penalty parameter would cause the loss function to become sharp, leading to inefficient training. 
All experiments are conducted on a single NVIDIA A$100$ GPU.

We remark that there are several potential enhancements for the aforementioned setup, such as using larger and deeper networks, training more epochs and so on. 
However, implementing these enhancements would require much more effort and resources. 
The current configuration is chosen by balancing efficiency and accuracy. 
Generally, it takes approximately $5$ minutes to discover a new eigenvalue in most cases.

\subsection{Fractional Laplacian in the $d$-dimensional unit ball}
We first calculate the eigenvalues of the fractional Laplace operator in $d$-dimensional unit balls, i.e., $\itOmega = B(0,1)$.
As we mentioned before, the exact eigenvalues of this problem are currently unknown. However, Dyda et al. show a method to calculate tight lower and upper bounds of the eigenvalues for this specific case\supercite{dyda2017eigenvalues}. We reproduce their method and obtain several bounds. By utilizing these bounds, we can confirm the first few digits of the exact eigenvalues. These inferred values are then used to test the accuracy of our method. 
In the following experiments, we calculate the relative errors by
\[
	e {:}= \frac{|\hat{\lambda} - \lambda^*|}{\lambda^*},
\]
where $\hat{\lambda}$ is our numerical result and $\lambda^*$ is the inferred value. 

We test our method for $d = 1,3,$ and $9$. For simplicity, we let the sampling domain $D = \itOmega$ and define the feature functions as
\begin{equation}
q_j(x){:}= \text{ReLU}\big(1 - \|x\|_2^2\big)^{p_j}.
\end{equation}
Here, $\text{ReLU}(x) = \max(0,x)$. These feature functions ensure that the output of the neural networks vanishes in $\itOmega^c$. The exponents $p_j$ are evenly spaced over the interval $[s,3]$ to capture both the sharp and the smooth behaviors near the boundary.

For the case $d=1$, the numerical results are summarized in Table 1. Our method successfully provides accurate values for the first $10$ eigenvalues, and their relative error are less than $0.2\%$.
For larger $s$, the results can be better, allowing us to calculate more eigenvalues while maintaining this level of precision.
Furthermore, we display the eigenfunctions in Figure 2 to demonstrate the singularity near the boundary. As $s$ decreases, the eigenfunctions become sharper near the boundary. 
\begin{table}
	\begin{center}
		{\small
			
			{\bf Table 1}\ \  Estimates of the eigenvalues of \eqref{problem_laplace} in $\itOmega = (-1,1)$.
			\vskip 1mm
			
			\begin{tabular}{ccllllll}
				\hline
				s&            & k=1      & k=2      & k=3      & k=4      & k=5      & k=10     \\ \hline
				& Exact      & 0.97259  & 1.09219  & 1.14732  & 1.18684  & 1.21655  & 1.31070  \\
				0.05 & Our        & 0.97261  & 1.09217  & 1.14735  & 1.18689  & 1.21665  & 1.31325  \\
				& Rel. error & 2.06\rme-05 & 1.83\rme-05 & 2.61\rme-05 & 4.38\rme-05 & 8.22\rme-05 & 1.95\rme-03 \\ \hline
				& Exact      & 0.97017  & 1.60154  & 2.02882  & 2.38716  & 2.69474  & 3.88845  \\
				0.25 & Our        & 0.97020  & 1.60148  & 2.02878  & 2.38761  & 2.69540  & 3.89149  \\
				& Rel. error & 3.09\rme-05 & 3.75\rme-05 & 1.97\rme-05 & 1.88\rme-04 & 2.45\rme-04 & 7.82\rme-04 \\ \hline
				& Exact      & 1.15777  & 2.75476  & 4.31680  & 5.89215  & 7.46018  & 15.3155  \\
				0.5  & Our        & 1.15780  & 2.75496  & 4.31666  & 5.89386  & 7.46028  & 15.3224  \\
				& Rel. error & 2.59\rme-05 & 7.26\rme-05 & 3.24\rme-05 & 2.90\rme-04 & 1.34\rme-05 & 4.51\rme-04 \\ \hline
				& Exact      & 1.59750  & 5.05976  & 9.59431  & 15.0188  & 21.1894  & 61.0924  \\
				0.75 & Our        & 1.59747  & 5.05971  & 9.59273  & 15.0225  & 21.1944  & 61.0977  \\
				& Rel. error & 1.88\rme-05 & 9.88\rme-06 & 1.65\rme-04 & 2.46\rme-04 & 2.36\rme-04 & 8.68\rme-05 \\ \hline
				& Exact      & 2.24406  & 8.59575  & 18.7168  & 32.4620  & 49.7200  & 186.450  \\
				0.95 & Our        & 2.24379  & 8.59504  & 18.7175  & 32.4626  & 49.7308  & 186.461  \\
				& Rel. error & 1.20\rme-04 & 8.26\rme-05 & 3.74\rme-05 & 1.94\rme-05 & 2.17\rme-04 & 5.90\rme-05 \\ \hline
			\end{tabular}
			
			\vskip 1mm
			\textit{Exact} show the first few digits of the exact eigenvalue.
		}
	\end{center}
	\label{Table1}
\end{table}
\begin{figure}
	\begin{center}
		
		\includegraphics[width=0.49\columnwidth]{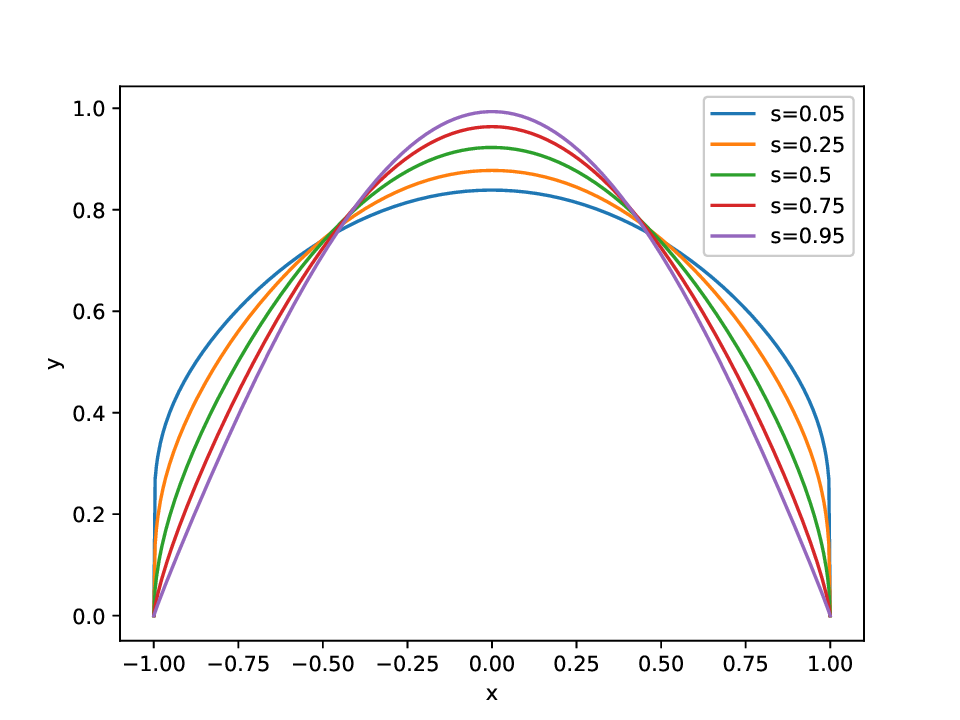}
		\includegraphics[width=0.49\columnwidth]{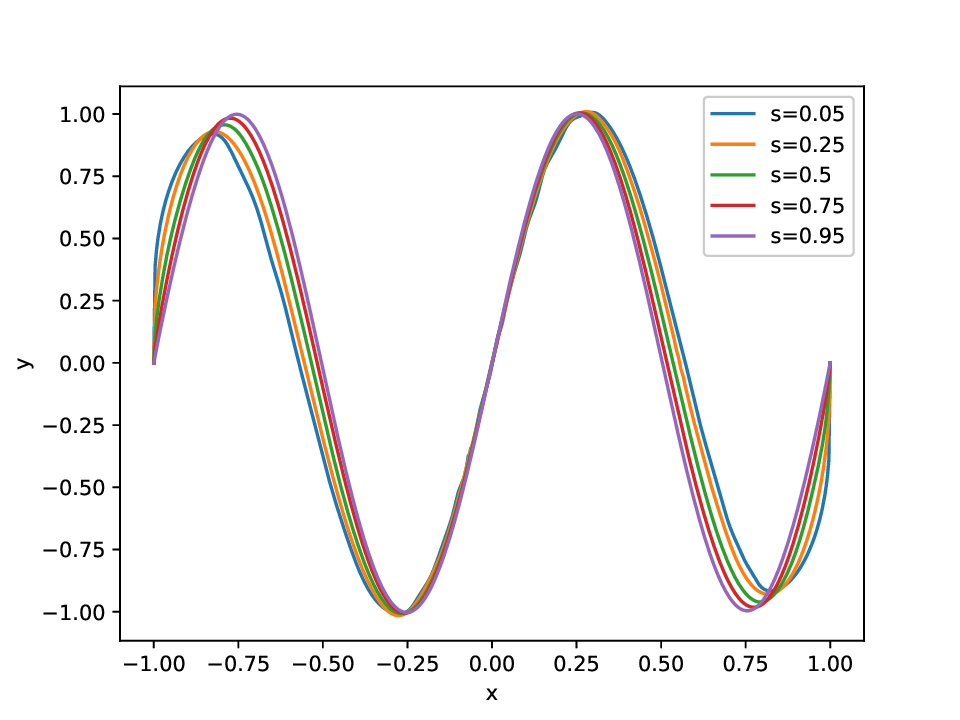}
		
		\small {\bf Figure 2}\ \ The first and the fourth eigenfunction of \eqref{problem_laplace} in $\itOmega = (-1,1)$.
		
	\end{center}
\end{figure}

When $d = 3$, we calculate the first $30$ eigenvalues. The relative error is less than $0.2\%$ for $s \geq 0.25$. 
However, the accuracy for $s = 0.05$ is not so satisfactory, compared to the higher order cases.
Their errors grow faster so that the results become unreliable earlier.
The numerical results are displayed in Table 2. 
We also calculate the result for $s=0.9999$.
Our numerical result is very close to the exact value while the eigenvalue of the fractional Laplacian converges to the eigenvalue of the common Laplacian as $s$ approaches $1$. 

In this example, some eigenvalues have a multiplicity greater than $1$. 
For these repeated eigenvalues, our method can identify all the mutually orthogonal eigenfunctions.
However, it is inefficient to find the next new eigenvalue with a different value if the multiplicity is too large. 
For example, in the case of the 9-dimensional ball, the smallest eigenvalue is simple. From the second to the tenth eigenvalue, they share the same value, while the next eigenvalue has a multiplicity of $44$. 
Consequently, it is very time-consuming for our method to find a new eigenvalue beyond these three eigenvalues.
The numerical results for the $9$-dimensional ball can be found in Table 3.

All these results indicate that solving the eigenvalue problem with a small fractional order is more challenging and the accuracy of the solution becomes worse earlier than the large fractional order cases. Additionally, it should be noted that the standard error of calculating \eqref{MC_eigen} caused by using the Monte Carlo method varies from $3\rme\text{-}05$ to $3\rme\text{-}04$ in all our experiments, depending on the complexity of the eigenmodes. Hence, the accuracy of these results has approached the limit of our method with the current configuration. 
\begin{table}
\begin{center}
{\small{\bf Table 2}\ \ Estimates of the eigenvalues of \eqref{problem_laplace} in the unit ball (d=3).
\vskip 1mm
\begin{tabular}{cclllllll}
\hline
s&            & k=1      & k=2      & k=3      & k=5      & k=10     & k=15     & k=30     \\ \hline
& Exact      & 1.092197 & 1.14300  & 1.14300  & 1.17687  & 1.18684  & 1.20274  & 1.23712  \\
0.05 & Our        & 1.092194 & 1.14303  & 1.14304  & 1.17757  & 1.18714  & 1.20518  & 1.26123*      \\
& Rel. error & 2.75\rme-06 & 2.62\rme-05 & 3.50\rme-05 & 5.95\rme-04 & 2.53\rme-04 & 2.03\rme-03 & 1.95\rme-02        \\ \hline
& Exact      & 1.601538 & 1.98571  & 1.98571  & 2.28647  & 2.38716  & 2.54207  & 2.92181  \\
0.25 & Our        & 1.601535 & 1.98569  & 1.98580  & 2.28750  & 2.38852  & 2.54400  & 2.92759  \\
& Rel. error & 1.87\rme-06 & 1.01\rme-05 & 4.53\rme-05 & 4.50\rme-04 & 5.70\rme-04 & 7.59\rme-04 & 1.98\rme-03 \\ \hline
& Exact      & 2.75476  & 4.12130  & 4.12130  & 5.40002  & 5.89215  & 6.63029  & 8.71829  \\
0.5  & Our        & 2.75498  & 4.12087  & 4.12114  & 5.40141  & 5.89405  & 6.63462  & 8.72610  \\
& Rel. error & 7.99\rme-05 & 1.04\rme-04 & 3.88\rme-05 & 2.57\rme-04 & 3.22\rme-04 & 6.53\rme-04 & 8.96\rme-04 \\ \hline
& Exact      & 5.05976  & 8.93319  & 8.93319  & 13.1781 & 15.0187 & 17.7566 & 26.5730  \\
0.75 & Our        & 5.06078  & 8.93035  & 8.93205  & 13.1780 & 15.0284 & 17.7742 & 26.5872  \\
& Rel. error & 2.02\rme-04 & 3.18\rme-04 & 1.28\rme-04 & 3.04\rme-06 & 6.40\rme-04 & 9.91\rme-04 & 5.33\rme-04 \\ \hline
& Exact      & 8.59575  & 17.0965  & 17.0965  & 27.5394 & 32.4619 & 39.8028 & 65.8034  \\
0.95 & Our        & 8.59548  & 17.0967  & 17.1011  & 27.5366 & 32.4666 & 39.8463 & 65.8256  \\
& Rel. error & 3.14\rme-05 & 1.05\rme-05 & 2.67\rme-04 & 1.01\rme-04 & 1.44\rme-04 & 1.09\rme-03 & 3.37\rme-04 \\ \hline

& Exact      & 9.86685  & 20.1840  & 20.1840  & 33.2050 & 39.4629 & 48.8111 & 82.6813  \\
0.9999 & Our        & 9.86767  & 20.1819  & 20.1830  & 33.2124 & 39.4798 & 48.8591 & 82.7274  \\
& Rel. error & 8.26\rme-05 & 1.02\rme-04 & 4.73\rme-05 & 2.23\rme-04 & 4.17\rme-04 & 9.82\rme-04 & 5.57\rme-04 \\ \hline

1 & Exact      & 9.86960  & 20.1907  & 20.1907  & 33.2175 & 39.4784 & 48.8312 & 82.7192  \\ \hline
\end{tabular}

\vskip 1mm
\textit{Exact} show the first few digits of the exact eigenvalue. * represents the error of this solution is too large and it should not be convinced.
}
\end{center}
\label{Table2}
\end{table}
\newpage
\begin{table}
\begin{center}
{\small
	
{\bf Table 3}\ \  Estimates of the eigenvalues of \eqref{problem_laplace} in the $9$-dimensional unit ball.

\begin{tabular}{cclllllll}
\hline
s &            & k=1      & k=2      & k=3      & k=5      & k=10     & k=11     & k=15     \\ \hline
& Exact      & 1.20274  & 1.22386  & 1.22386  & 1.22386  & 1.22386  & 1.24179  & 1.24179  \\
0.05 & Our        & 1.20275  & 1.22392  & 1.22393  & 1.22395  & 1.22475  & 1.24361  & 1.24411  \\
& Rel. error & 8.31\rme-06 & 4.90\rme-05 & 5.72\rme-05 & 7.35\rme-05 & 7.27\rme-04 & 1.47\rme-03 & 1.87\rme-03 \\ \hline
& Exact      & 2.54207  & 2.76833  & 2.76833  & 2.76833  & 2.76833  & 2.97357  & 2.97357  \\
0.25 & Our        & 2.54212  & 2.76834  & 2.76855  & 2.76884  & 2.78037  & 2.98222  & 2.98411  \\
& Rel. error & 1.97\rme-05 & 3.61\rme-06 & 7.95\rme-05 & 1.84\rme-04 & 4.35\rme-03 & 2.91\rme-03 & 3.54\rme-03 \\ \hline
& Exact      & 6.63029  & 7.82911  & 7.82911  & 7.82911  & 7.82911  & 9.00556  & 9.00556  \\
0.5  & Our        & 6.62929  & 7.82897  & 7.82986  & 7.83152  & 7.87850  & 9.02421  & 9.03107  \\
& Rel. error & 1.51\rme-04 & 1.79\rme-05 & 9.58\rme-05 & 3.08\rme-04 & 6.31\rme-03 & 2.07\rme-03 & 2.83\rme-03 \\ \hline
& Exact      & 17.7566  & 22.6391  & 22.6391  & 22.6391  & 22.6391  & 27.8025  & 27.8025  \\
0.75 & Our        & 17.7617  & 22.6390  & 22.6414  & 22.6581  & 22.7577  & 27.8484  & 27.9111  \\
& Rel. error & 2.87\rme-04 & 4.42\rme-06 & 9.81\rme-05 & 8.39\rme-04 & 5.24\rme-03 & 1.65\rme-03 & 3.91\rme-03 \\ \hline
& Exact      & 39.8028  & 53.8038  & 53.8039  & 53.8038  & 53.8038  & 69.4807  & 69.4807  \\
0.95 & Our        & 39.8146  & 53.8138  & 53.8427  & 53.8439  & 54.0359  & 69.6889  & 69.7573  \\
& Rel. error & 2.96\rme-04 & 1.86\rme-04 & 7.22\rme-04 & 7.45\rme-04 & 4.31\rme-03 & 3.00\rme-03 & 3.98\rme-03 \\ \hline
\end{tabular}

\vskip 1mm
\textit{Exact} show the first few digits of the exact eigenvalue.

}
\end{center}
\label{Table3}
\end{table}

\subsection{Fractional Schrödinger operator}
Next, we solve the eigenvalue problem of the fractional Schrödinger operator with a potential function $V(x)$. To demonstrate the effectiveness of our method, we conduct two tests in one-dimensional intervals and solve a three-dimensional problem with an inverse square potential in a unit ball.

The problem domains for the first two examples are $\itOmega = (-1,1)$. We also let $D = \itOmega$ and define the feature functions as
\begin{equation}
q_j(x){:}= \text{ReLU}\big(1 - x_1^2\big)^{p_j}.
\end{equation}
The exponents $p_j$ are also evenly spaced over the interval $[s,3]$. 
In the first example, the potential function we employed is $V(x) = x^2/2$. 
We test our method for different $s$ and the numerical results are presented in Table 4.
By comparing our results with those obtained from \cite{bao2020jacobi}, we observe that they are very close to each other.
The relative differences are less than $0.06\%$ for the first $10$ eigenvalues in all circumstances.

The second example involves a distinct potential function $V(x) = 50x^2 + \sin(2 \pi x)$. 
We calculate the first few eigenvalues and show our estimates in Table 5.
Additionally, we plot the first six eigenfunctions in Figure 3, revealing that the shapes and singularity of the eigenfunctions of the fractional Schrödinger operator differ from those of the fractional Laplacian. It is clear that some eigenfunctions do not exhibit singularity near the boundary with this potential. But, our method still successfully identifies them.
\newpage
\begin{table}
	\begin{center}
		{\small{\bf Table 4}\ \  Estimates of the eigenvalues of \eqref{problem_schrodinger} with $V(x) = x^2/2$ in $\itOmega = (-1,1)$.
			\vskip 1mm
			
			\begin{tabular}{llllllll}
				\hline
				s    &       & k=1      & k=2      & k=3      & k=4      & k=5      & k=10     \\ \hline
				& Our   & 1.05992  & 1.76847  & 2.19033  & 2.55183  & 2.85805  & 4.05264  \\
				0.25 & Ref. & 1.05995  & 1.76850  & 2.19047  & 2.55226  & 2.85848  & 4.05477  \\
				& Diff. & 2.83E-05 & 1.70E-05 & 6.39E-05 & 1.69E-04 & 1.50E-04 & 5.26E-04 \\ \hline
				& Our   & 1.24024  & 2.91807  & 4.48137  & 6.05866  & 7.62650  & 15.4822  \\
				0.5  & Ref. & 1.24036  & 2.91792  & 4.48124  & 6.05836  & 7.62828  & 15.4813  \\
				& Diff. & 9.68E-05 & 5.14E-05 & 2.90E-05 & 4.95E-05 & 2.33E-04 & 5.81E-05 \\ \hline
				& Our   & 1.67073  & 5.21206  & 9.75501  & 15.1826  & 21.3543  & 61.2587  \\
				0.75 & Ref. & 1.67054  & 5.21184  & 9.75495  & 15.1818  & 21.3573  & 61.2629  \\
				& Diff. & 1.14E-04 & 4.22E-05 & 6.15E-06 & 5.27E-05 & 1.40E-04 & 6.86E-05 \\ \hline
				& Our   & 2.31064  & 8.73900  & 18.8735  & 32.6231  & 49.8832  & 186.616 \\
				0.95 & Ref. & 2.31063  & 8.73878  & 18.8749  & 32.6228  & 49.8845  & 186.667 \\
				& Diff. & 4.33E-06 & 2.52E-05 & 7.42E-05 & 9.20E-06 & 2.61E-05 & 2.74E-04 \\ \hline
			\end{tabular}
			
			\vskip 1mm
			\textit{Ref.} represents the reference values given by \cite{bao2020jacobi}.
			\textit{Diff.} represents the relative difference between these two results.
			
		}
	\end{center}
	\label{Table4}
\end{table}
\begin{table}
	\begin{center}
		{\small			
			{\bf Table 5}\ \  Estimates of the eigenvalues of \eqref{problem_schrodinger} with $V(x) = 50x^2+\sin(2\pi x)$ in $\itOmega = (-1,1)$.
			\vskip 1mm
			
			\begin{tabular}{lllllll}
				\hline
				s    & k=1     & k=2     & k=3     & k=4     & k=5     & k=10    \\ \hline
				0.25 & 2.17977 & 3.90875 & 4.74208 & 5.50170 & 6.06804 & -*       \\
				0.5  & 3.67234 & 8.61490 & 11.9706 & 15.0614 & 17.7531 & 29.2724 \\
				0.75 & 5.31594 & 14.5127 & 22.3722 & 30.0024 & 37.4224 & 78.2848 \\
				0.95 & 6.71889 & 20.0334 & 33.6763 & 48.9012 & 66.5992 & 203.348 \\ \hline
			\end{tabular}
			\vskip 1mm
			* represents it fails to generate a solution due to the accumulation of the compuataional error.
		}
	\end{center}
	\label{Table5}
\end{table}
\begin{figure}
	\begin{center}
		\centering
		\includegraphics[scale=0.28]{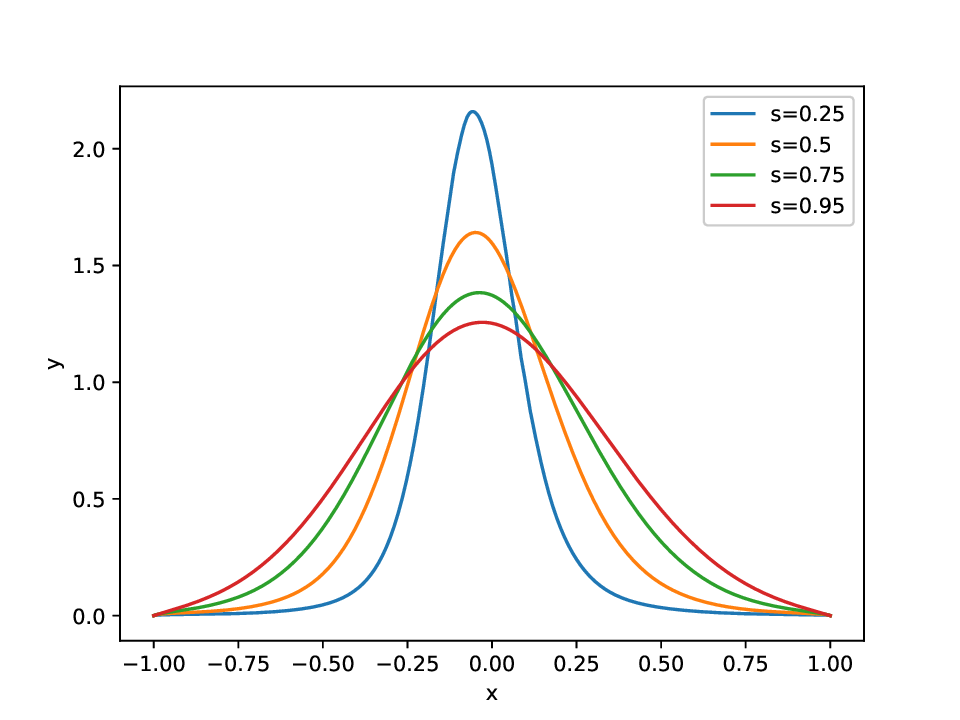}
		\includegraphics[scale=0.28]{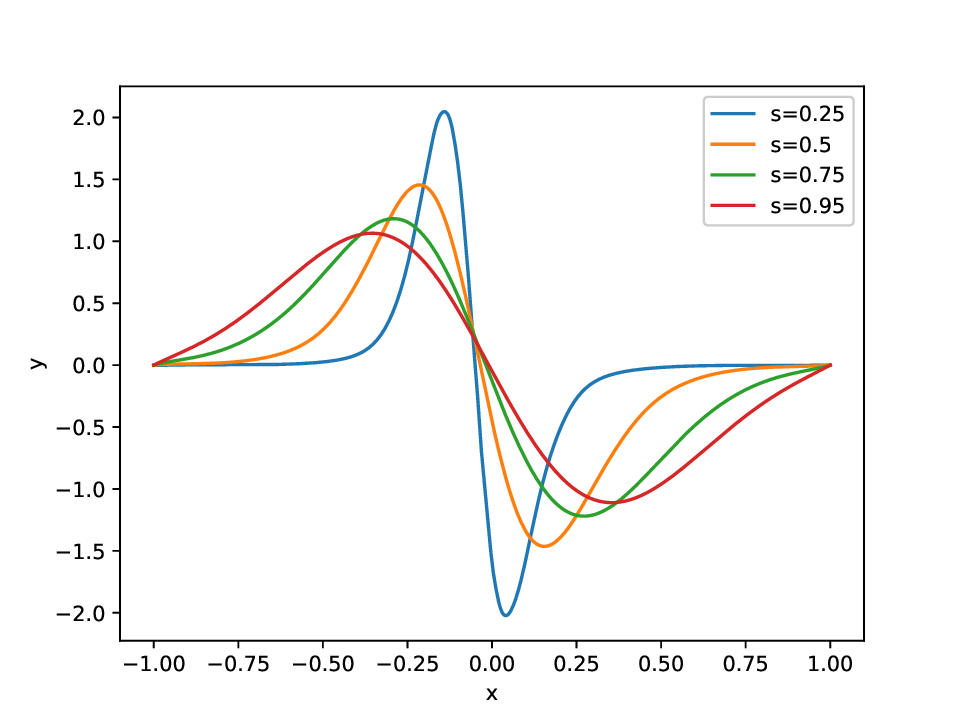}
		\includegraphics[scale=0.28]{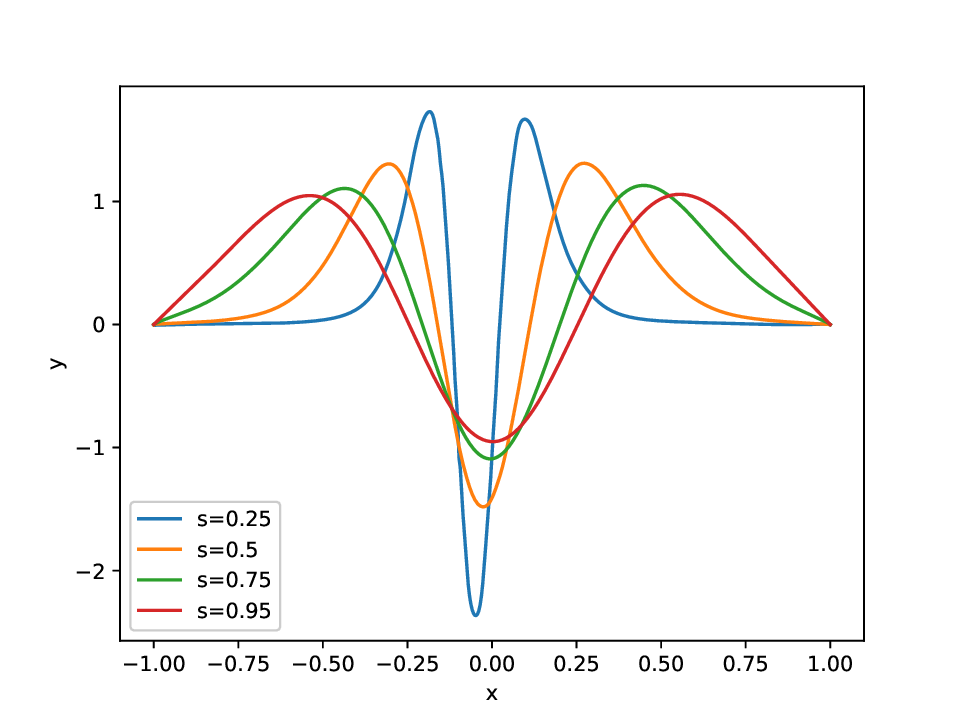}
		
		\includegraphics[scale=0.28]{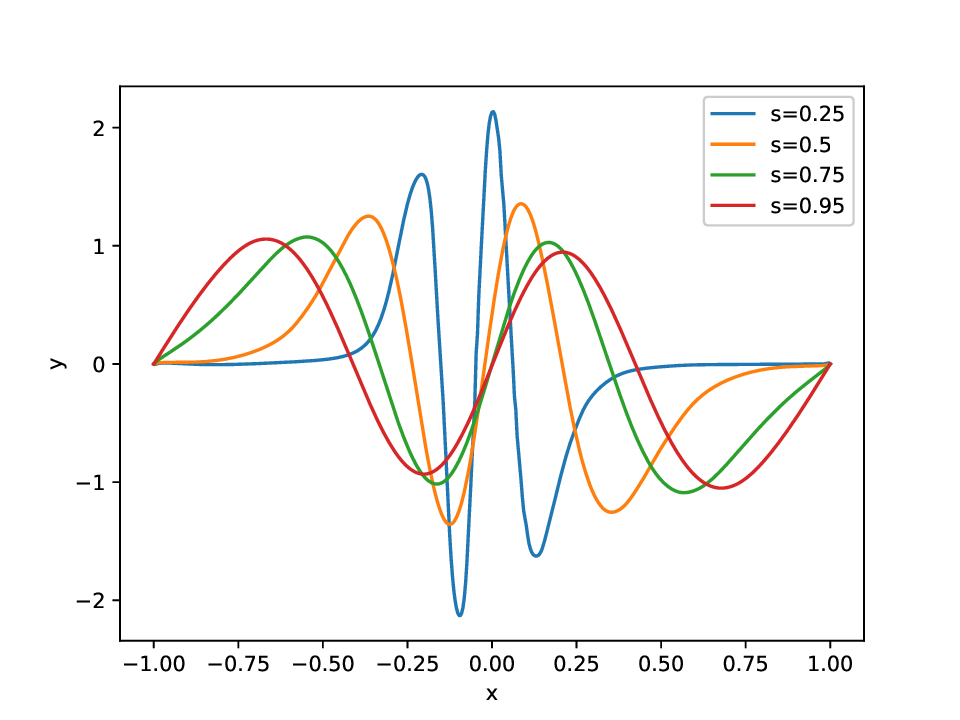}
		\includegraphics[scale=0.28]{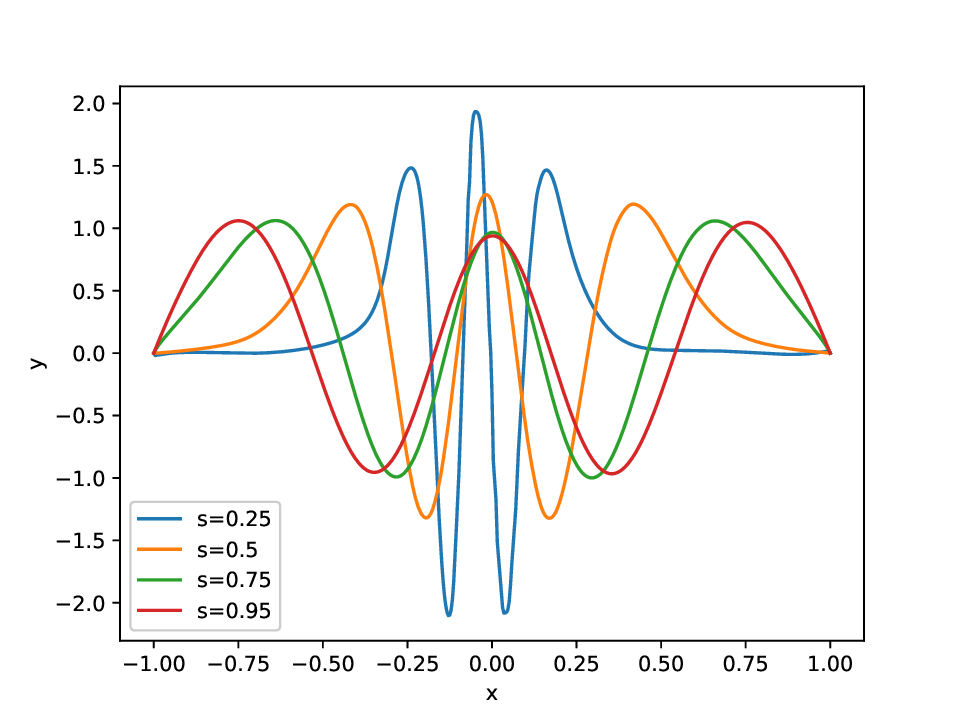}
		\includegraphics[scale=0.28]{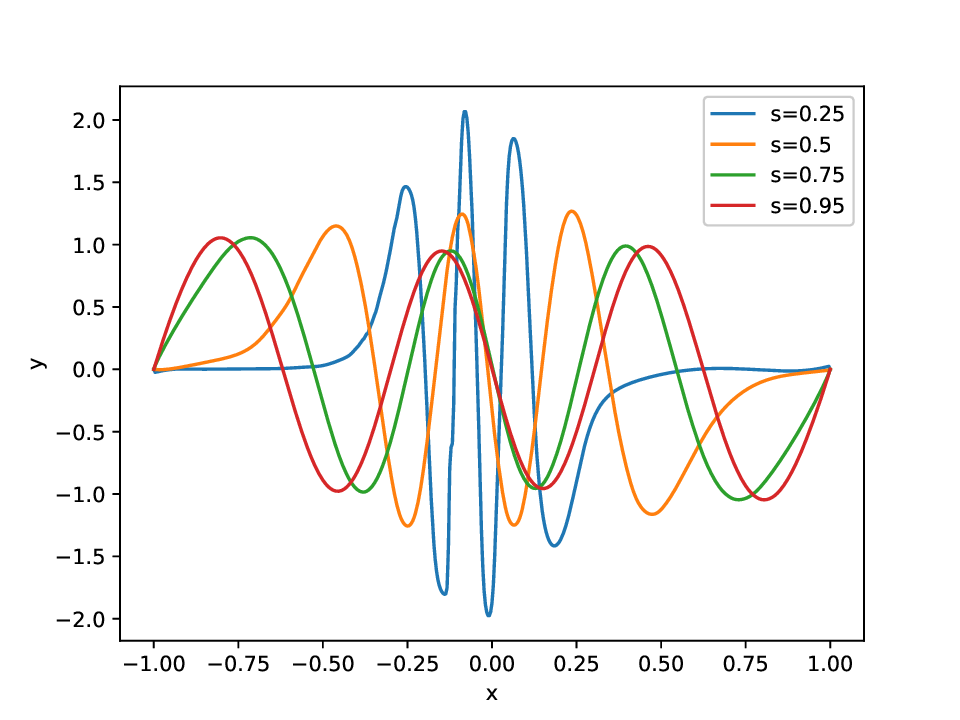}
		
		{\small 
			{\bf Figure 3}\ \ The first six eigenfunctions of \eqref{problem_schrodinger} with $V(x) = 50x^2+\sin(2\pi x)$ in $\itOmega = (-1,1)$
		}
	\end{center}
\end{figure}

Last, we solve the problem \eqref{problem_schrodinger} in the unit ball with an inverse square potential 
\[
V(x) = \frac{1}{2(x_1^2+x_2^2+x_3^2)}.
\]
The feature functions we used are the same as those in $V(x) = 0$. 
Figure 4 shows the eigenvalues with these two different potential functions.
In all other cases, the order of the eigenvalues is independent of $s$, i.e., they would not exchange for different $s$. However, in this case, the order of the eigenvalues changes. 
As $s$ decreases, the value of the first $7$-fold eigenvalue is no longer greater than the value of the second single eigenvalue, and the value of the first $9$-fold eigenvalue is no longer greater than the value of the second triple eigenvalue.
\begin{figure}
\begin{center}

\subfigure[potential $V(x) = \frac{1}{2(x_1^2+x_2^2+x_3^2)}$]
{\includegraphics[scale=0.44]{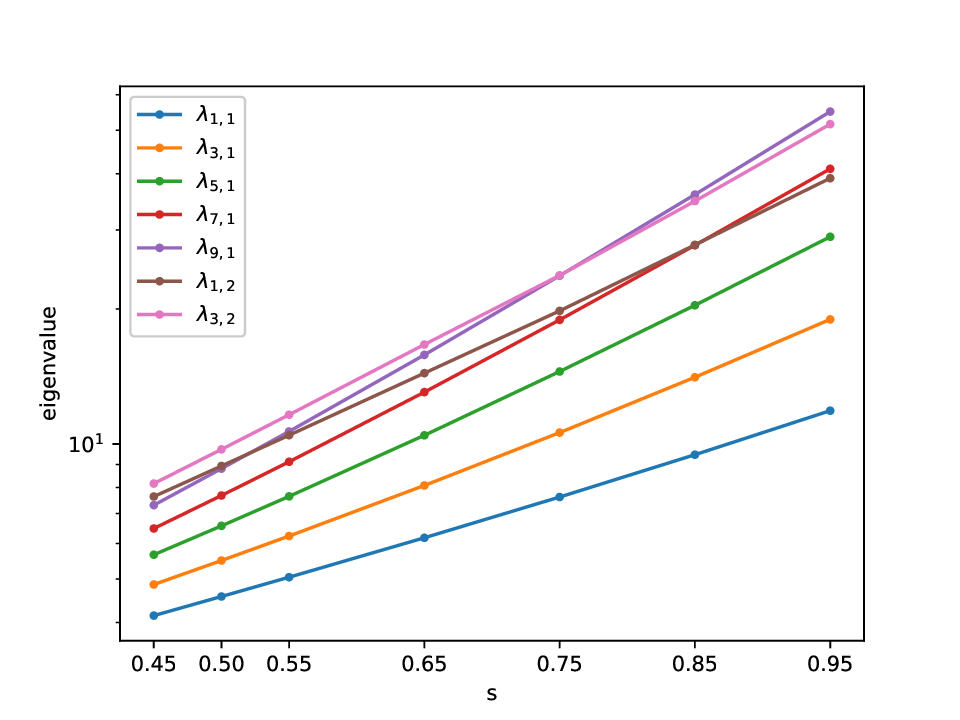}}
\subfigure[potential $V(x) = 0$]{\includegraphics[scale=0.44]
	{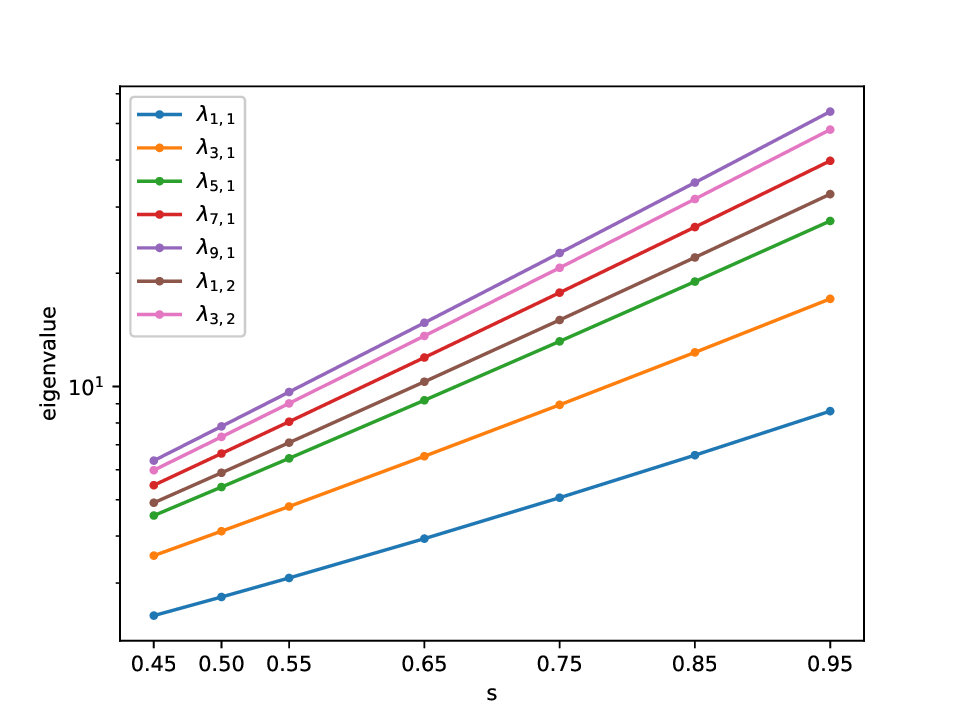}}

{\small 
	\parbox[t]{1.7cm}{\bf Figure 4}
	\parbox[t]{10.5cm}{Eigenvalues of \eqref{problem_schrodinger} in the $3$-dimensional unit ball with different potentials,
	$\lambda_{a,1}$ represents the value of the first eigenvalue which has a multiplicity of $a$ 
	while $\lambda_{a,2}$ represents the value of the second eigenvalue which has a multiplicity of $a$.}
}

\end{center}
\end{figure}
\subsection{Fractional Laplacian in general domains}
In this subsection, we focus on calculating the eigenvalues of the fractional Laplace operator over general domains. 
To validate our method, we compare the results with those computed by the finite element method in \cite{FEM_borthagaray2018finite}.

Firstly, we consider the problem in $\itOmega = [-1,1]^2$ and let the sampling domain $D = \itOmega$. The feature function is defined as
\begin{equation}
q_j(x){:}= \bigg[ \text{ReLU}\big(1-x_1^2\big) \cdot \text{ReLU}\big(1-x_2^2\big) \bigg]^{p_j}.
\end{equation}
Here, $p_j$ are also evenly spaced over the interval $[s,3]$.
We calculate the first eigenvalue for different $s$ and the outcomes are summarized in Table 6.
The results demonstrate that our method outperforms the finite element method over the finest grid, and these values are very close to the extrapolated values obtained through Richardson extrapolation using the finite element method's results.
It is worth noting that our approach enables us to obtain the corresponding eigenfunctions, which are not provided by the extrapolation method. 
Table 6 also demonstrates the efficiency of our method, as it only takes around $5$ minutes to calculate a new eigenmode.
We further compute more eigenvalues, and the results are presented in Table 7.
The multiplicity of each computed eigenvalue is the same as that of the Laplacian, which is well-known (see, for example, \cite{Laplacian_eigenvalue_square3_liu2013verified}).
\begin{table}
\begin{center}
{\small
	
{\bf Table 6}\ \  Estimates of the first eigenvalue of \eqref{problem_laplace} in the square $[-1,1]^2$.
\vskip 1mm

\begin{tabular}{cllllc}
\hline
s    & Our     & Extrapolated & FEM    &  & Our Time(s) \\ \cline{2-4} \cline{6-6} 
0.05 & 1.04054 & 1.0405      & 1.0412 &  &  290.2      \\
0.25 & 1.28129 & 1.2813      & 1.2844 &  &  292.4      \\
0.5  & 1.83440 & 1.8344      & 1.8395 &  &  291.8      \\
0.75 & 2.88721 & 2.8872      & 2.8921 &  &  290.6      \\
0.95 & 4.40568 & 4.4062      & 4.4083 &  &  296.5      \\ \hline
\end{tabular}

\vskip 1mm
\textit{Extrapolated} indicates the extrapolated values with FEM results in \cite{FEM_borthagaray2018finite}, \textit{FEM} represents the eigenvalues calculated by the Finite Element Method over the finest mesh in \cite{FEM_borthagaray2018finite}.

}
\end{center}
\label{Table6}
\end{table}

\begin{table}
	\begin{center}
		{\small
			
			{\bf Table 7}\ \  Estimates of the eigenvalues of \eqref{problem_laplace} in the square $[-1,1]^2$.
			\vskip 1mm
			
			\begin{tabular}{cccccccc}
				\hline
				s    & k=1     & k=2,3   & k=4     & k=5,6   & k=7,8   & k=9,10  & k=11    \\ \hline
				0.05 & 1.04054 & 1.10942 & 1.14281 & 1.15823 & 1.17429 & 1.19215 & 1.19685 \\
				0.25 & 1.28129 & 1.72109 & 1.97902 & 2.09892 & 2.26684 & 2.43361 & 2.48065 \\
				0.5  & 1.83440 & 3.14066 & 4.08501 & 4.59306 & 5.30757 & 6.09787 & 6.31421 \\
				0.75 & 2.88721 & 6.06243 & 8.79700 & 10.4447 & 12.8430 & 15.7654 & 16.5448 \\
				0.95 & 4.40568 & 10.6589 & 16.7414 & 20.7257 & 26.6512 & 34.4497 & 36.4295 \\ \hline
			\end{tabular}
			
		}
	\end{center}
	\label{Table7}
\end{table}

Next, we turn to the problem in an $L$-shaped domain $\itOmega = [-1,1]^2 \backslash [0,1]^2$. Since $\itOmega$ is not convex, we let the sampling domain $D = [-1,1]^2$. We employ two types of feature functions. The first-type function is 
\begin{equation}
q_{1,j}(x){:}= \max\bigg\{ 
\text{ReLU}\big(-x_1(x_1 + 1)\big) \text{ReLU}\big(1-x_1^2\big), \ 
\text{ReLU}\big(-x_2(x_2 + 1)\big) \text{ReLU}\big(1-x_2^2\big) \bigg\}^{p_j}.
\end{equation} 
Similar to the previous case, the exponents $p_j$ are evenly spaced over the interval $[s,3]$. It is well known that the solution of Laplace's equation over the L-shaped domain exhibits a singularity of type $r^{2/3}f(\theta)$ at the corner. We suspect that certain eigenfunctions may also display a corner singularity. To capture this singularity, we use another type of feature function:
\begin{equation}
q_{2,j}(x){:}= B\big(2r(x)\big) 
\ \sin\bigg(\frac{2}{3} \ \text{ReLU}\big(\theta(x) - \frac{\pi}{2}\big)\bigg) \ r(x)^{t_j}.
\end{equation} 
Here, $r(x)$ represents the distance between point $x$ and the corner, while the angle $\theta(x)$ is defined as the angle between the positive $x$-axis and the line connecting the point $x$ and the corner in a counterclockwise direction.
The exponents $t_j$ are evenly spaced over the interval $[2/3, 3/2]$.
$B(\cdot)$ is a bump function defined as
\begin{equation}
B(x) = 
\left \{
\begin{aligned}
&\exp\bigg(-\frac{1}{1-x^2}\bigg), & x\in (-1,1), \\
&0, &\text{otherwise}.
\end{aligned}
\right .
\end{equation}

In this example, we let the width of the network be $60$ and compare two different numerical schemes to demonstrate that incorporating the knowledge about the singularity near the corner can further enhance accuracy for $s > 2/3$. The first scheme we used, denoted as Scheme $A$ in the following paragraphs, utilizes $40$ first-type feature functions and $20$ second-type feature functions, while the second scheme, referred to as Scheme $B$, only employs the first-type feature functions. All other settings for these two schemes remain the same.

According to Table 8, Scheme $A$ provides lower estimates than Scheme $B$ when $s=0.7$ and $s=0.9$. This can be attributed to the fact that Scheme $A$ incorporates more knowledge about the singularity. However, the results of these two schemes are similar when $s \leq 0.6$, and both schemes outperform the results in \cite{FEM_borthagaray2018finite}.

We plot some eigenfunctions of Scheme $A$ in Figure 5, which indicates that the eigenfunctions exhibit similar shapes. But, the eigenfunctions with smaller fractional orders change more sharply near the boundary. Consequently, the absolute value of these functions tends to approach zero in a more narrow region. 
These behaviors are consistent with the previous example of the fractional Laplacian. 
In Table 9, we present our estimates for the first few eigenvalues using Scheme $A$. It shows that the multiplicities of these eigenvalues are consistent with those of the Laplacian\supercite{Laplacian_eigenvalue_Lshaped2_yuan2009bounds}.
\begin{table}
\begin{center}
{\small

{\bf Table 8}\ \ Estimates of the first eigenvalue of \eqref{problem_laplace} in the $L$-shaped domain $[-1,1]^2 \backslash [0,1]^2$.
\vskip 1mm

\begin{tabular}{cllllccc}
		\hline
s   & Our - A & Our - B & Extrapolate & FEM    &  & Time(s) - A & Time(s) - B \\ \cline{2-5} \cline{7-8} 
0.1 & 1.14145 & 1.14145 & 1.1413      & 1.1434 &  & 397.7       & 319.9           \\
0.3 & 1.59621 & 1.59609 & 1.5956      & 1.6025 &  & 398.1       & 322.1           \\
0.5 & 2.43299 & 2.43316 & 2.4322      & 2.4440 &  & 399.1       & 323.5           \\
0.6 & 3.09453 & 3.09478 & 3.0936      & 3.1072 &  & 397.4       & 325.4           \\
0.7 & 4.00864 & 4.00952 & 4.0069      & 4.0228 &  & 398.1       & 324.1           \\
0.9 & 7.08512 & 7.09517 & 7.0790      & 7.0975 &  & 399.5       & 322.3           \\ \hline
\end{tabular}

\vskip 1mm
\textit{Extrapolated} indicates the extrapolated values with FEM results in \cite{FEM_borthagaray2018finite}, \textit{FEM} represents the eigenvalues calculated by the Finite Element Method over the finest mesh in \cite{FEM_borthagaray2018finite}.
}
\end{center}
\label{Table8}
\end{table}

\begin{table}
\begin{center}
{\small

{\bf Table 9}\ \  Estimates of the eigenvalues of \eqref{problem_laplace} in the $L$-shaped domain $[-1,1]^2 \backslash [0,1]^2$.
\vskip 1mm

\begin{tabular}{ccccccccc}
\hline
s    & k=1     & k=2     & k=3     & k=4     & k=5     & k=6     & k=7     & k=8,9   \\ \hline
0.05 & 1.06508 & 1.11541 & 1.13495 & 1.16620 & 1.16646 & 1.18297 & 1.18921 & 1.20039 \\
0.25 & 1.45545 & 1.77523 & 1.92973 & 2.19286 & 2.20084 & 2.35759 & 2.41845 & 2.52001 \\
0.5  & 2.43299 & 3.38167 & 3.95538 & 5.00859 & 5.10676 & 5.84385 & 6.12008 & 6.56155 \\
0.75 & 4.59315 & 6.91928 & 8.59204 & 11.9203 & 12.4673 & 15.2133 & 16.2296 & 17.6734 \\
0.95 & 8.25174 & 12.9116 & 16.6466 & 24.5501 & 26.3578 & 33.8279 & 36.5487 & 40.0902 \\ \hline
\end{tabular}

}
\end{center}
\label{Table9}
\end{table}

\newpage
\begin{figure}
\begin{center}
\centering
\includegraphics[scale=0.21]{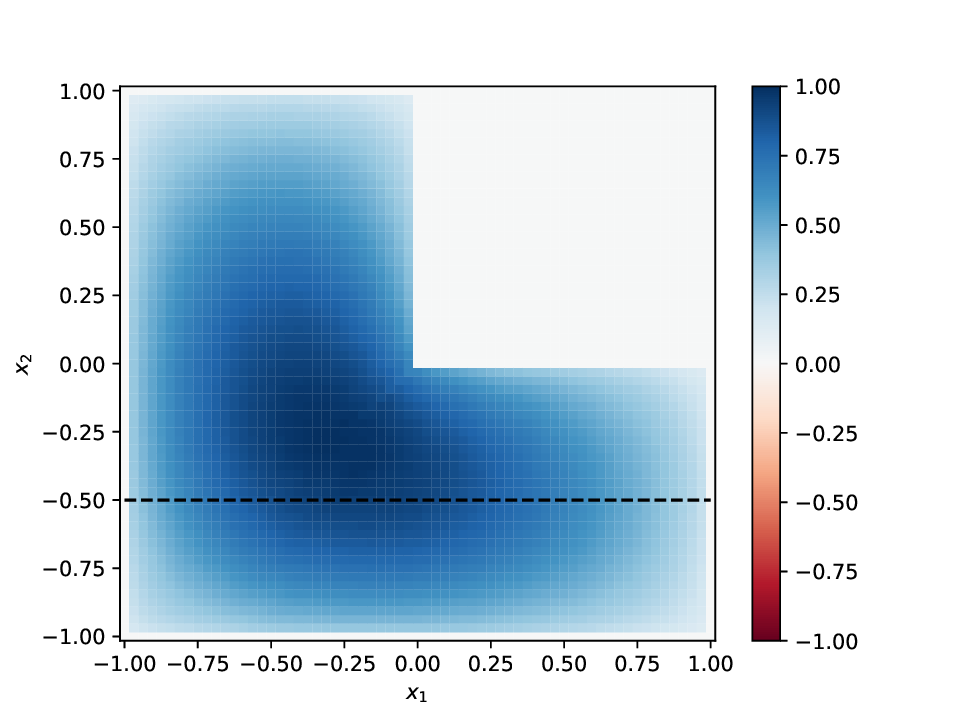}
\includegraphics[scale=0.21]{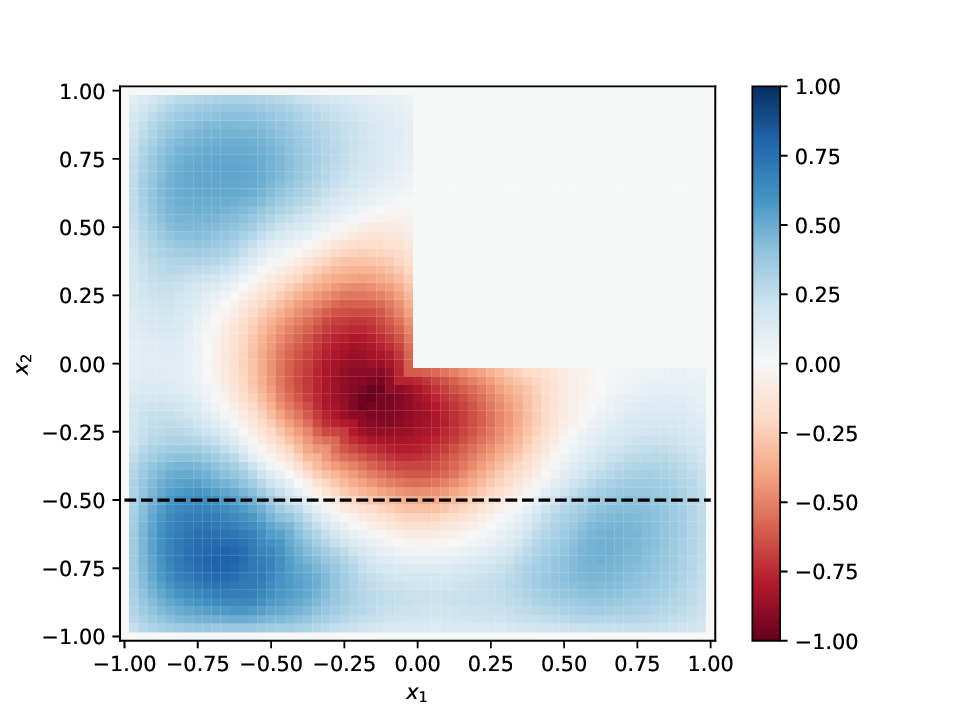}
\includegraphics[scale=0.21]{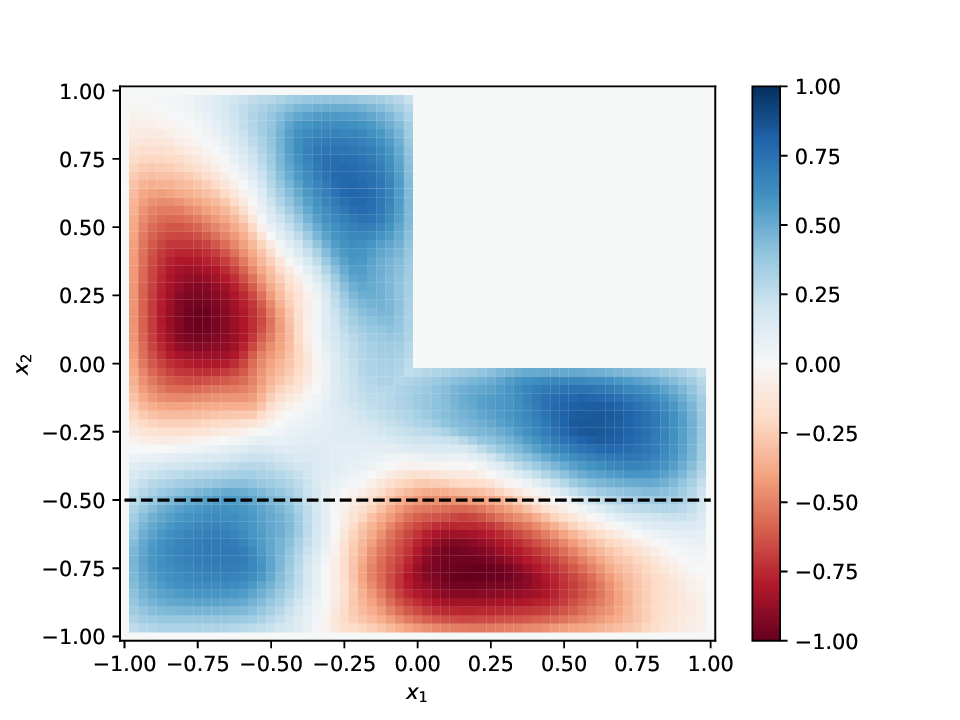}
\includegraphics[scale=0.21]{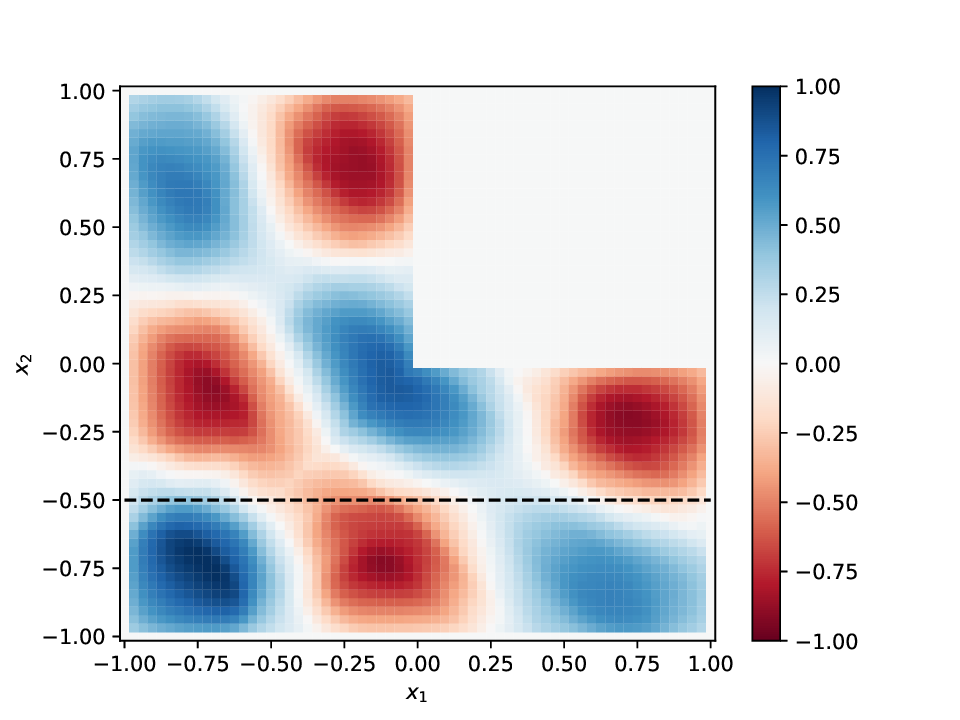}

\includegraphics[scale=0.21]{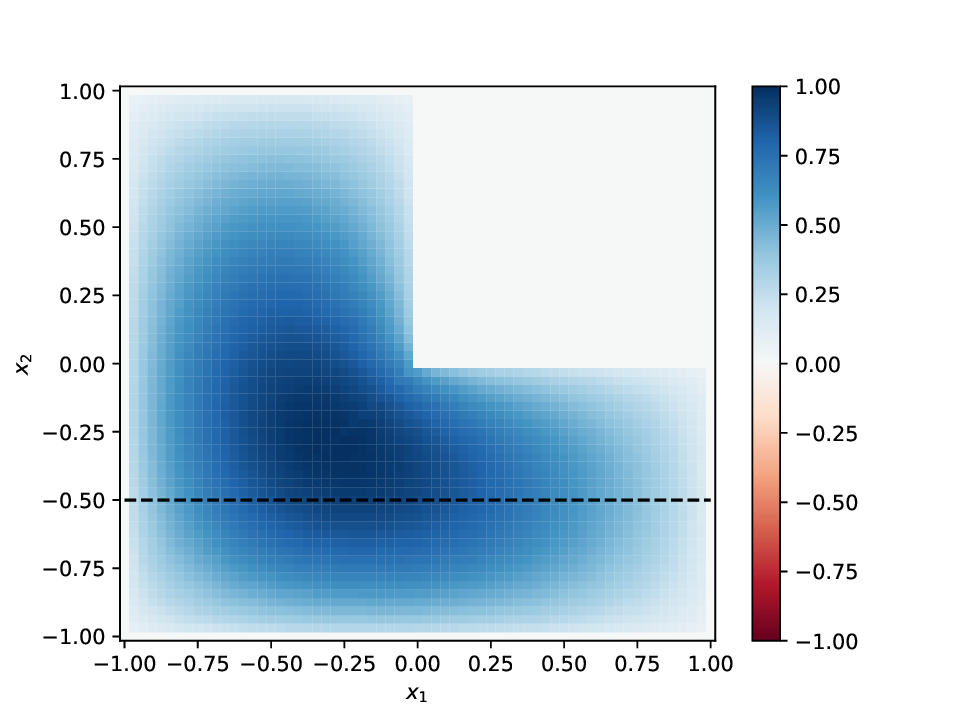}
\includegraphics[scale=0.21]{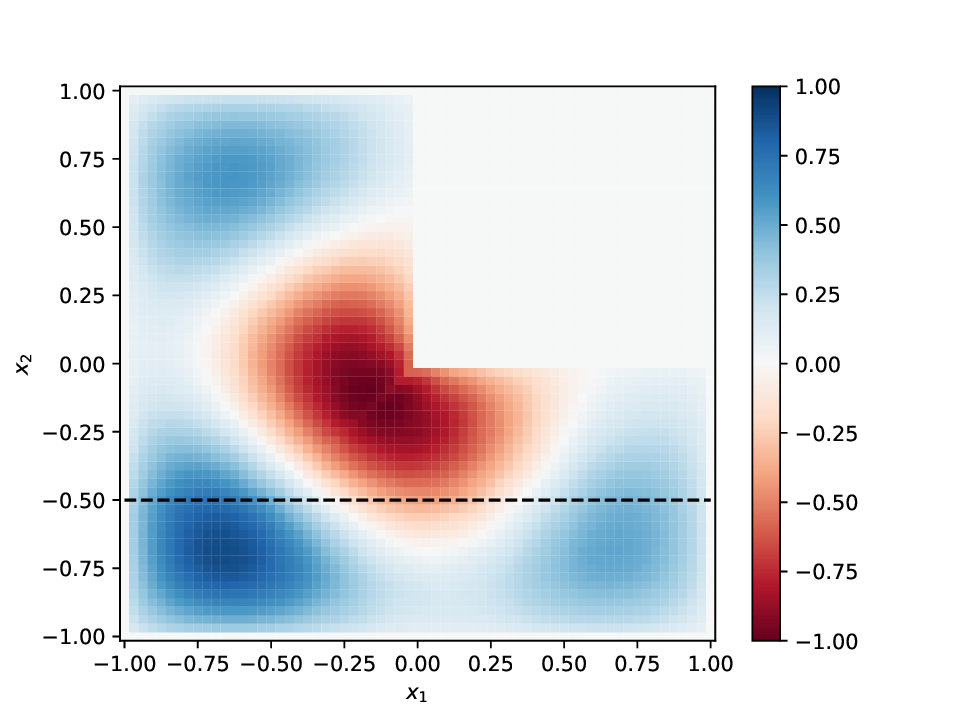}
\includegraphics[scale=0.21]{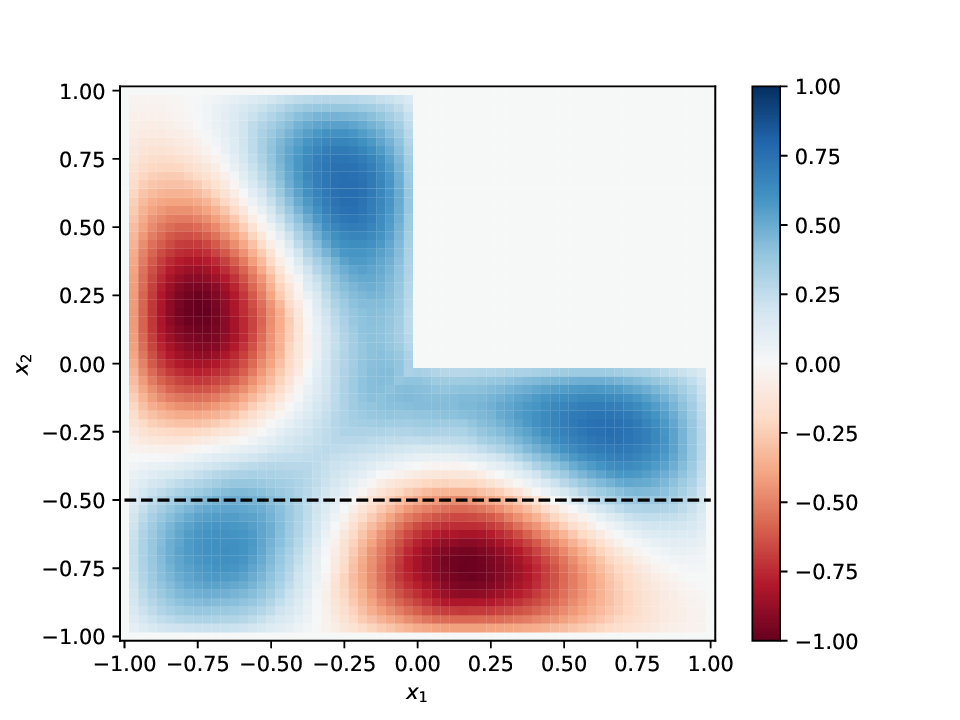}
\includegraphics[scale=0.21]{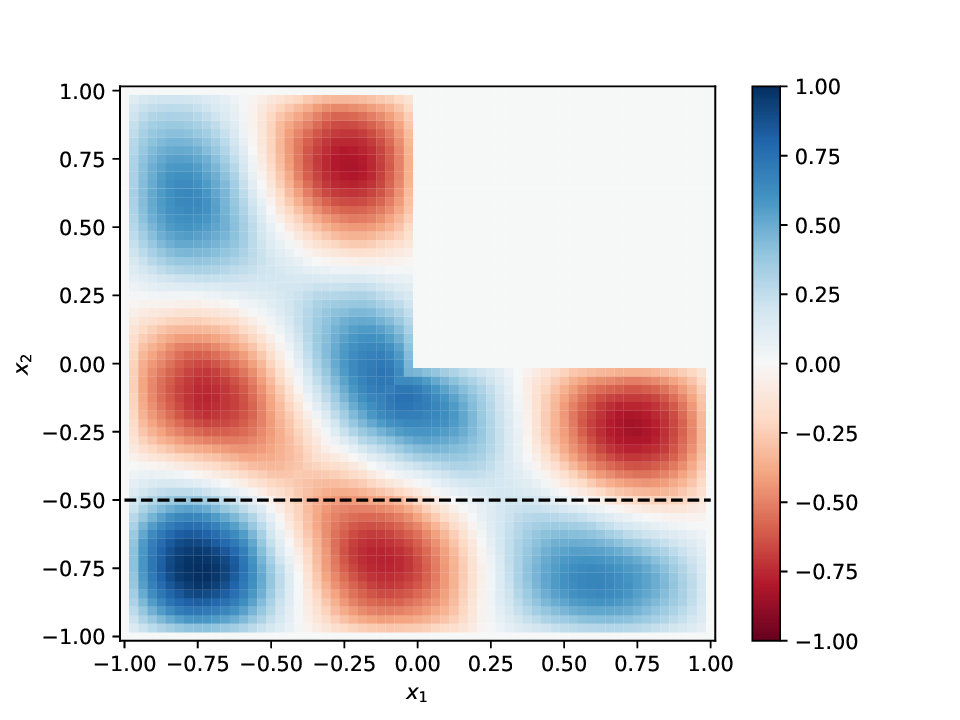}

\includegraphics[scale=0.21]{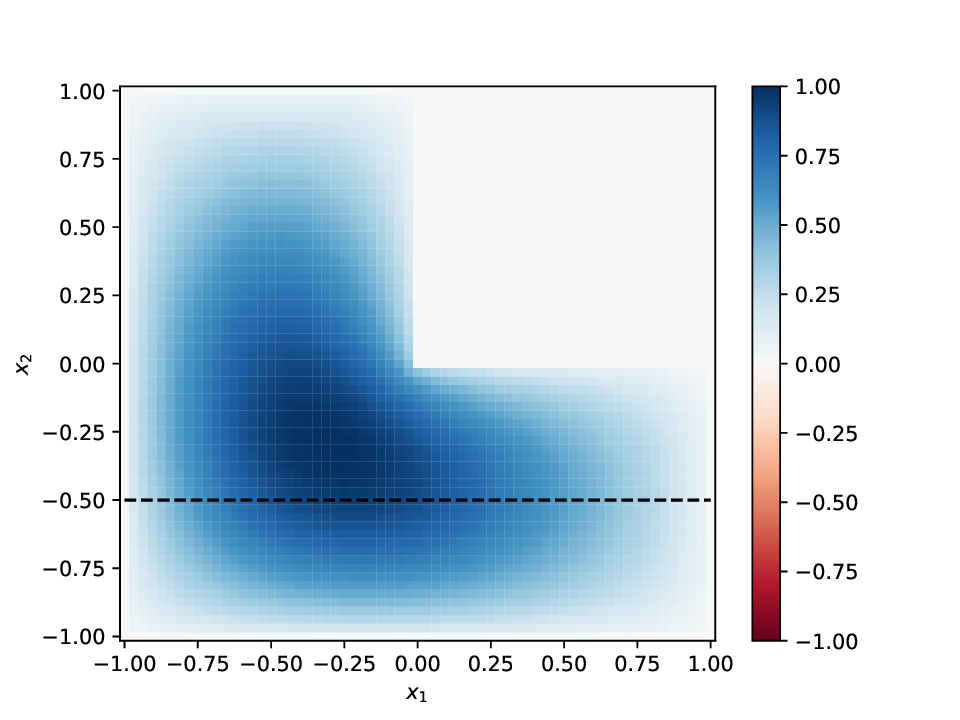}
\includegraphics[scale=0.21]{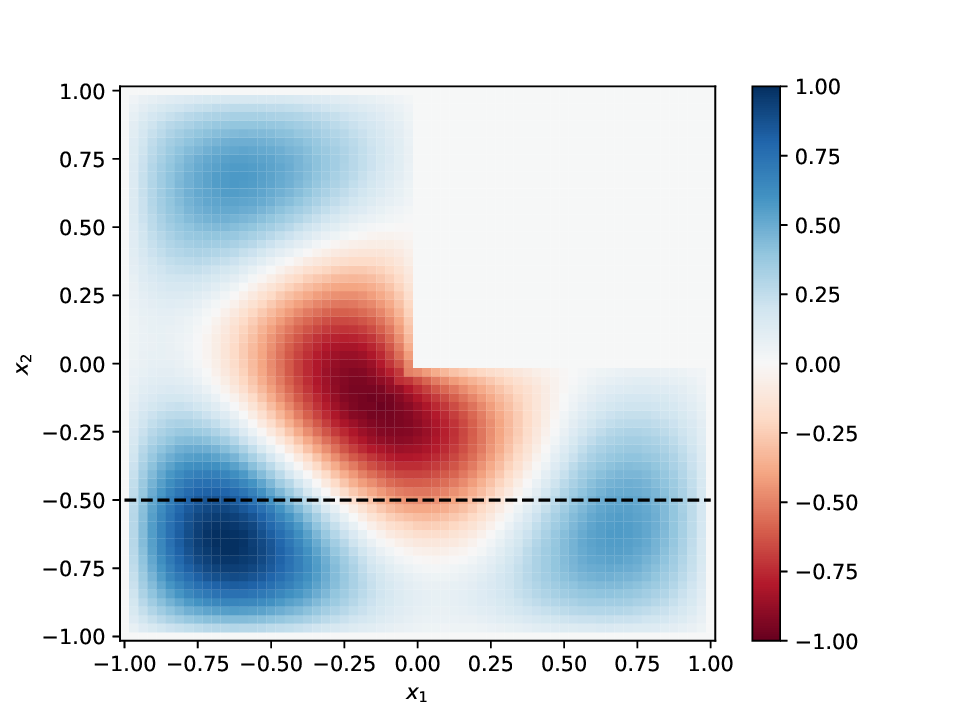}
\includegraphics[scale=0.21]{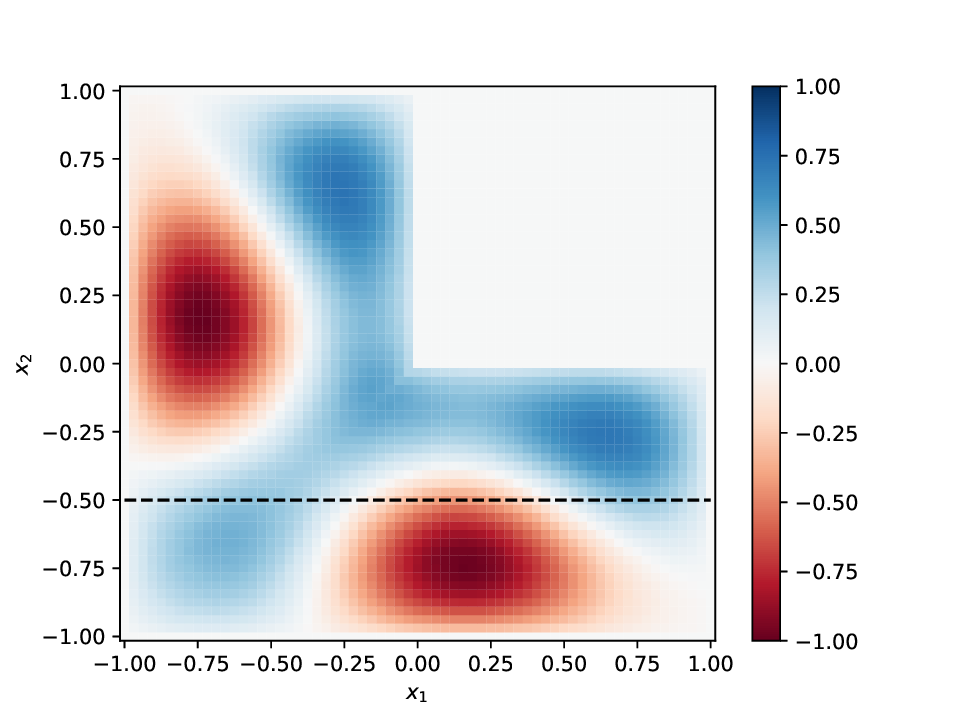}
\includegraphics[scale=0.21]{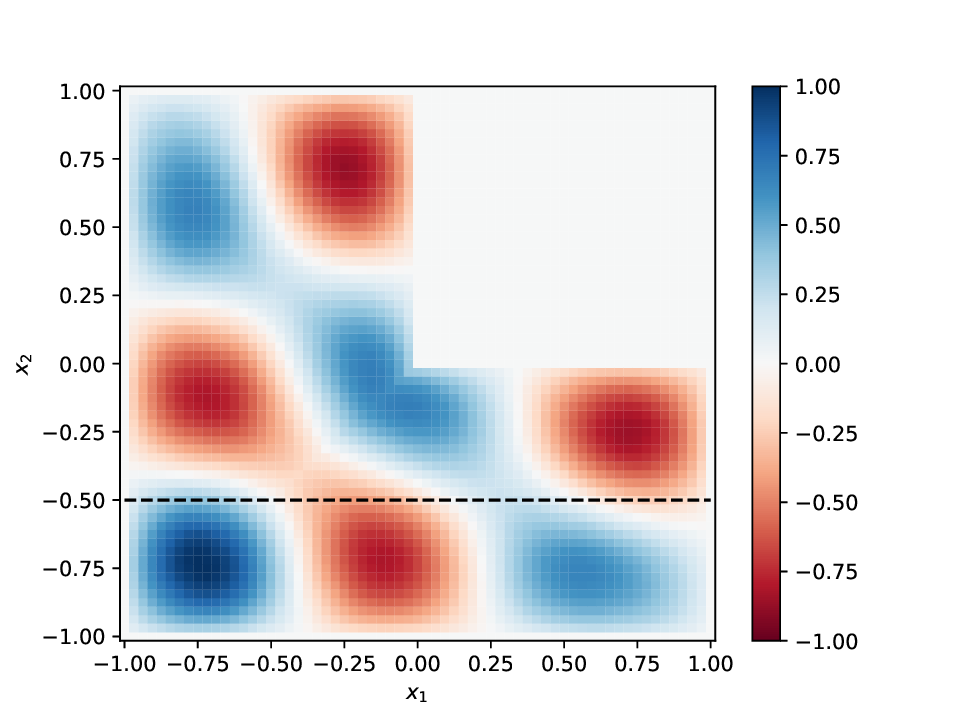}

\includegraphics[scale=0.21]{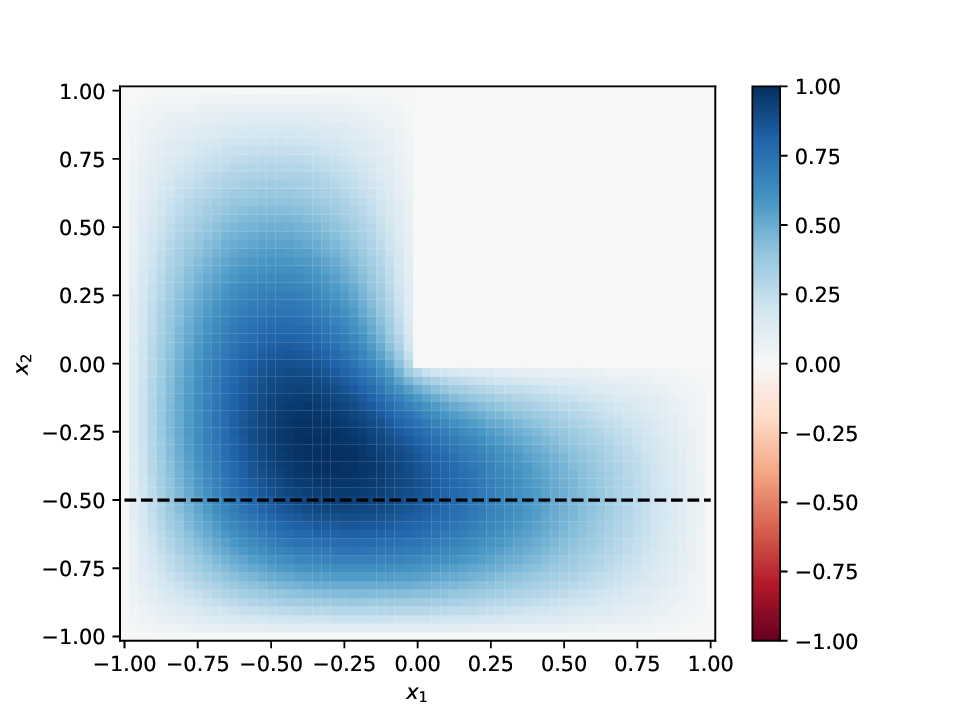}
\includegraphics[scale=0.21]{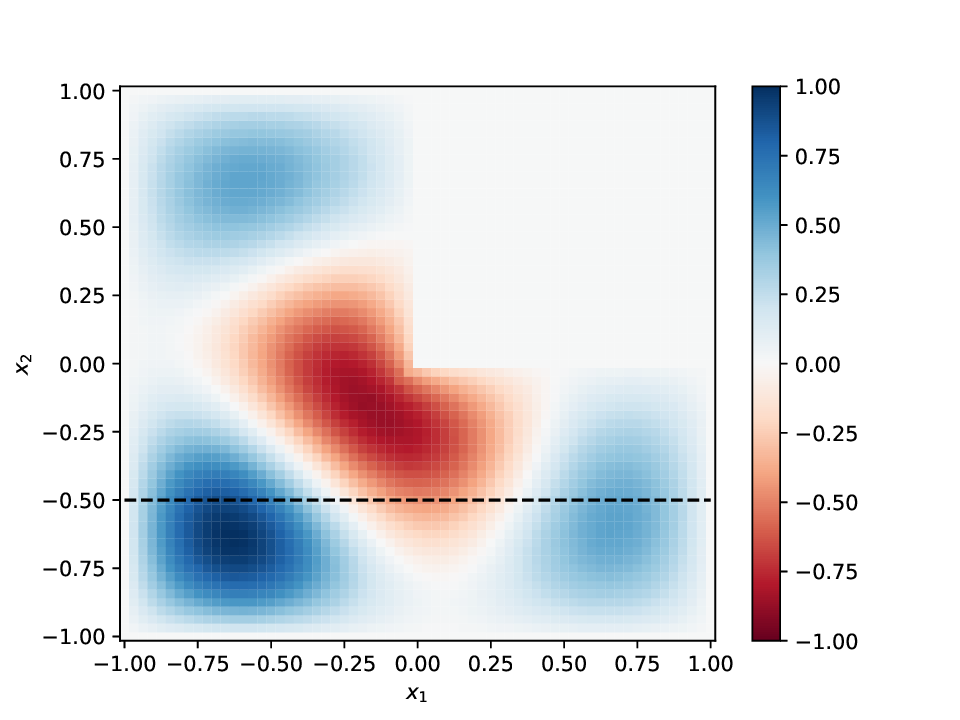}
\includegraphics[scale=0.21]{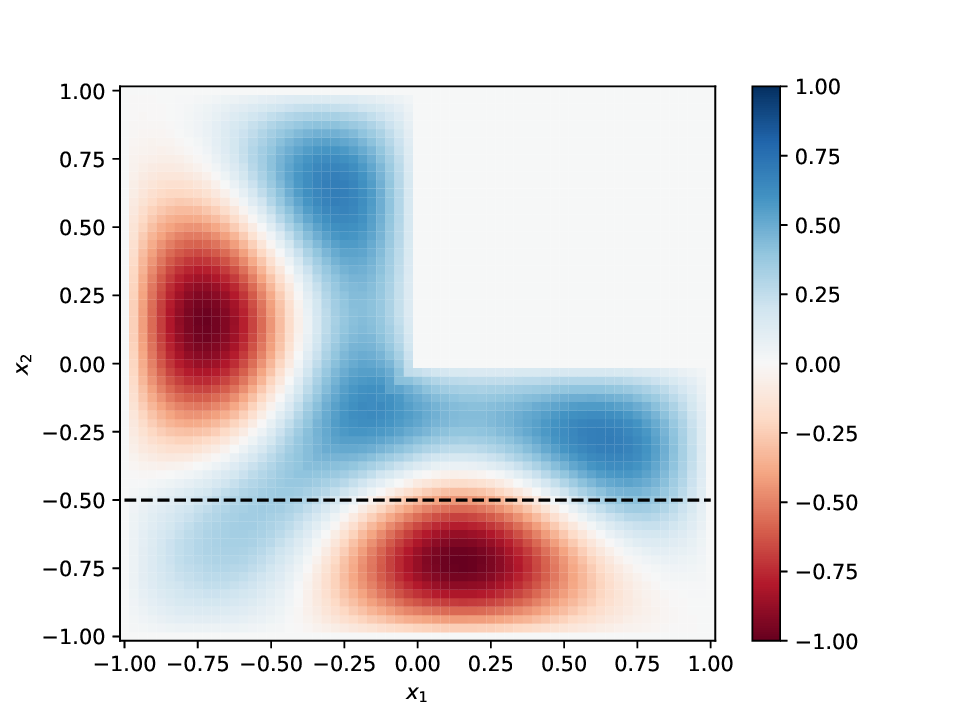}
\includegraphics[scale=0.21]{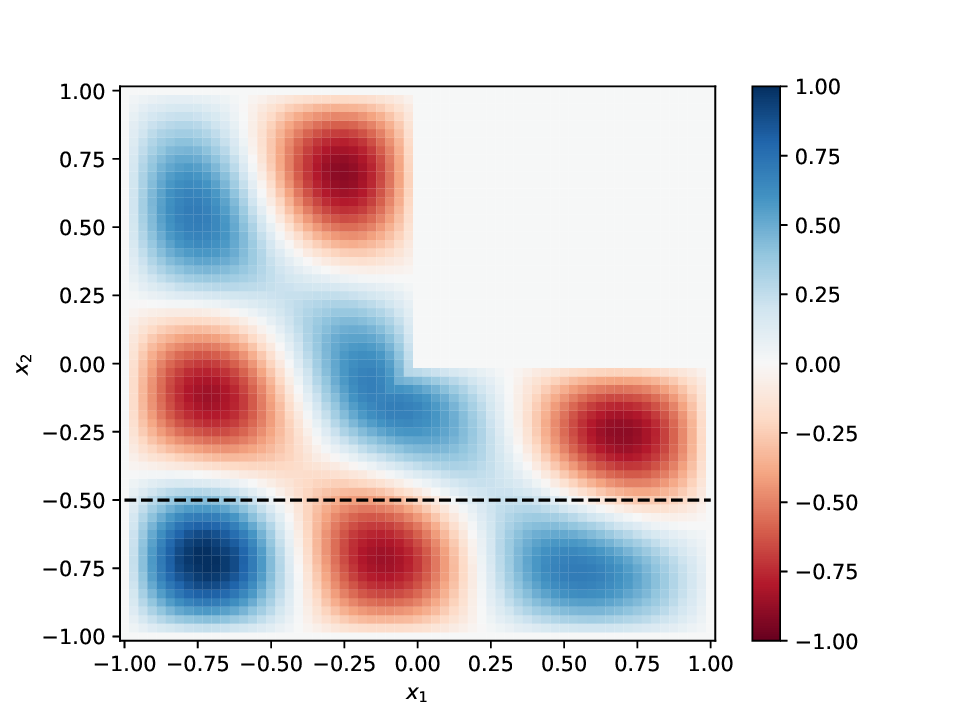}

\includegraphics[scale=0.21]{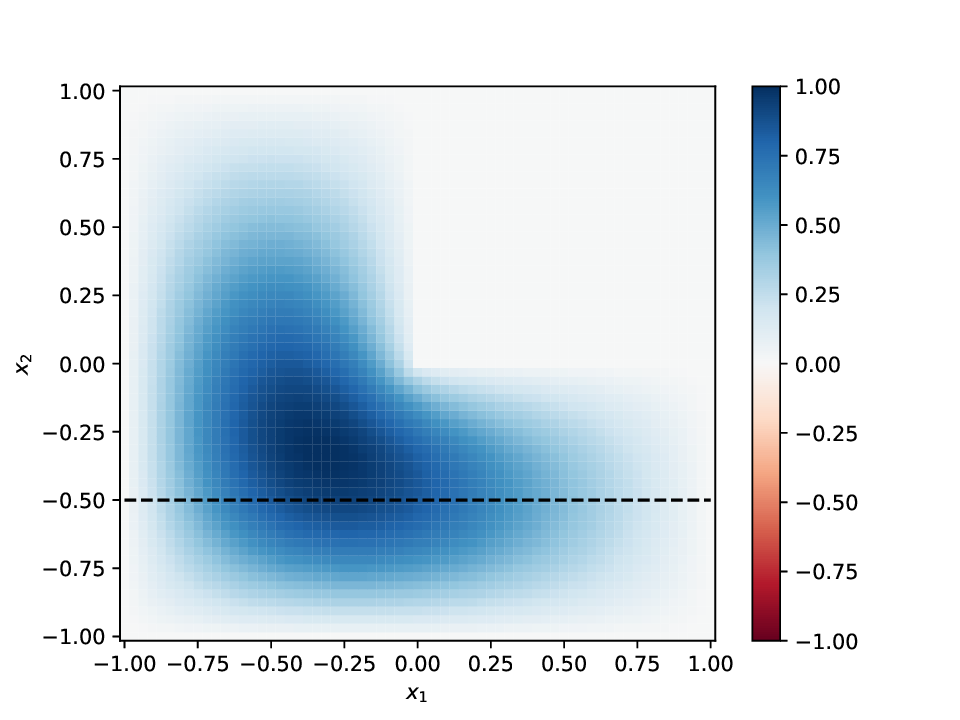}
\includegraphics[scale=0.21]{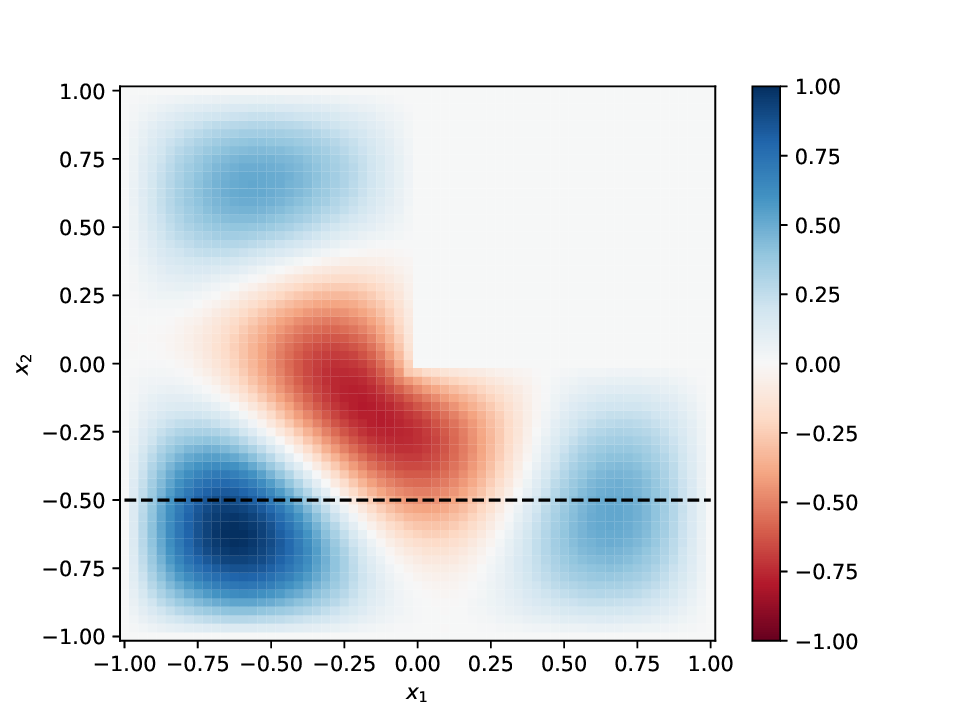}
\includegraphics[scale=0.21]{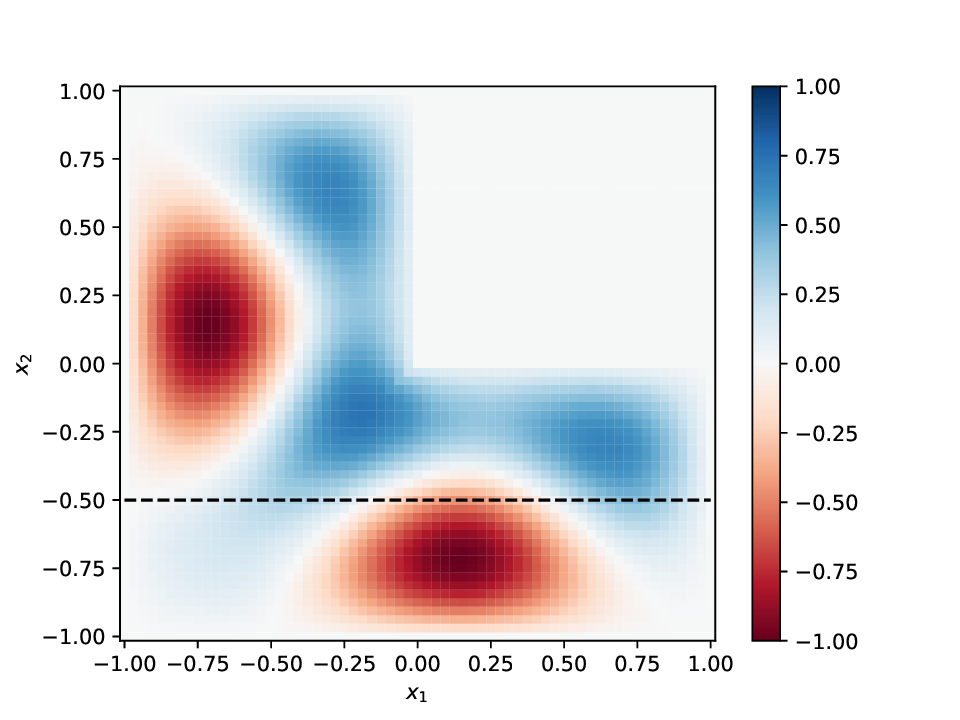}
\includegraphics[scale=0.21]{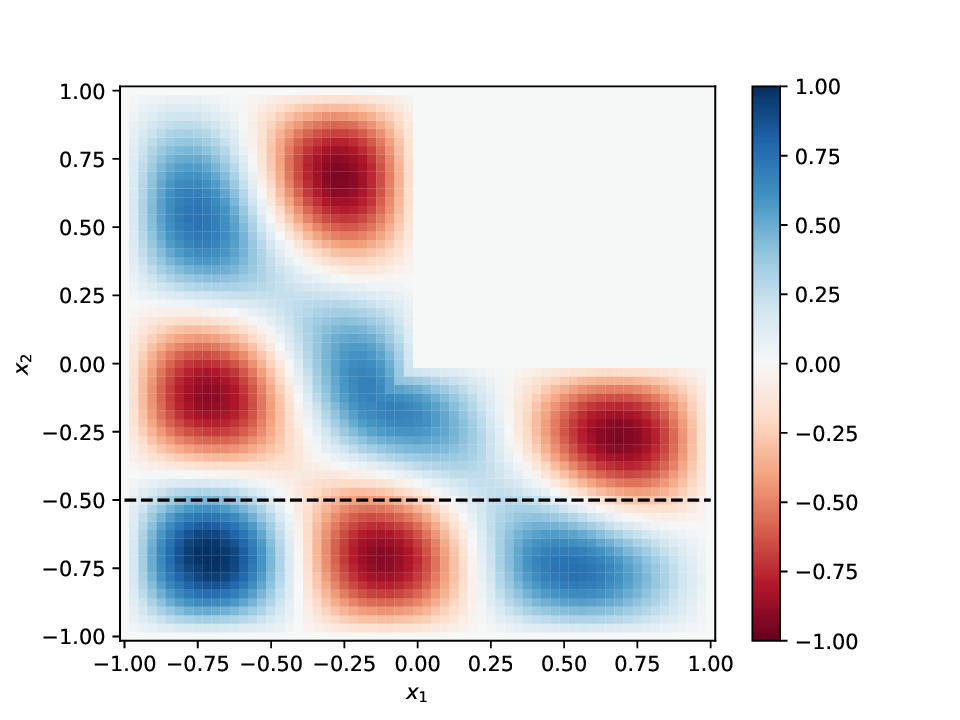}

\includegraphics[scale=0.21]{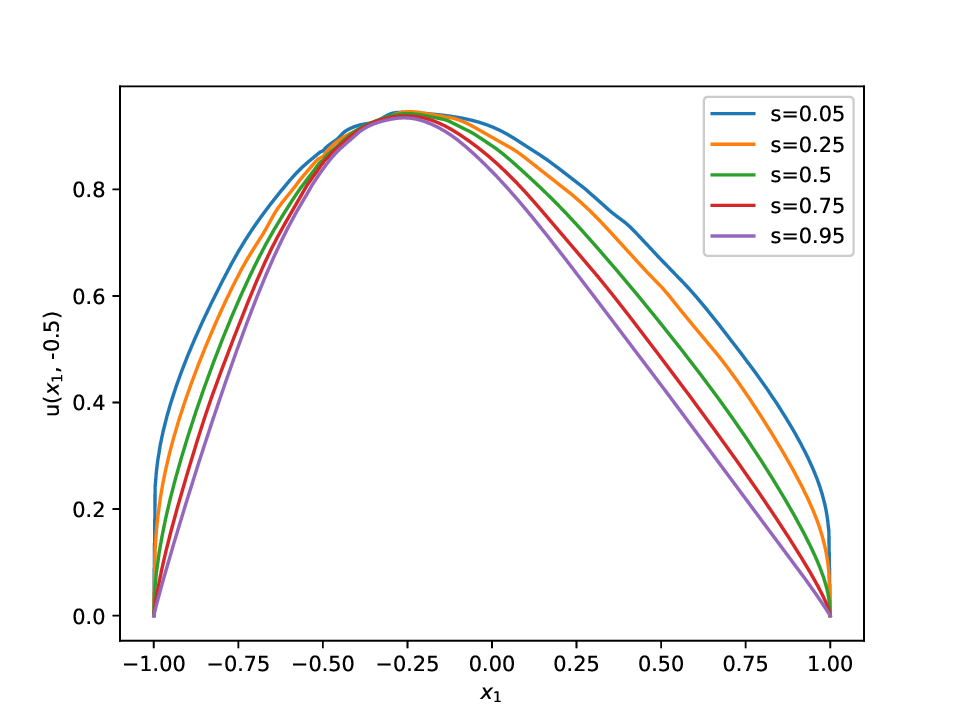}
\includegraphics[scale=0.21]{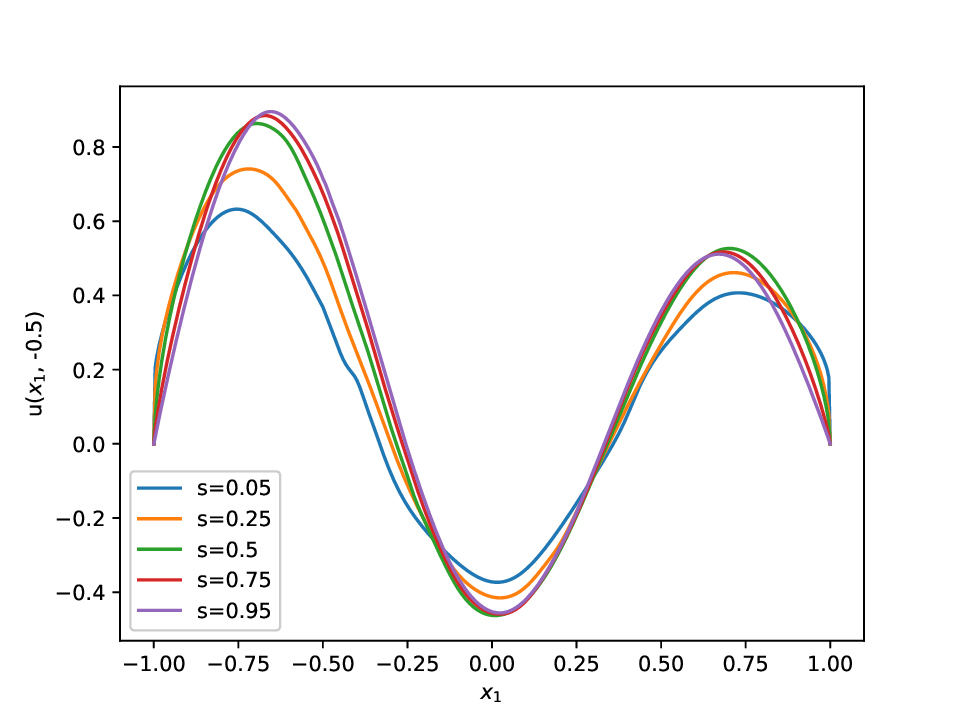}
\includegraphics[scale=0.21]{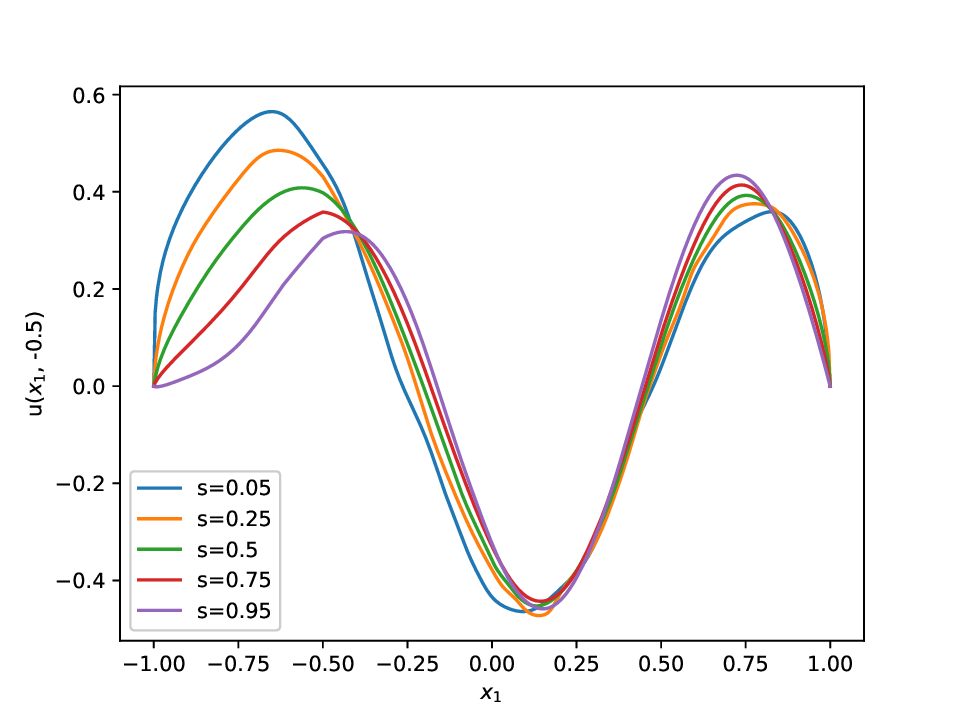}
\includegraphics[scale=0.21]{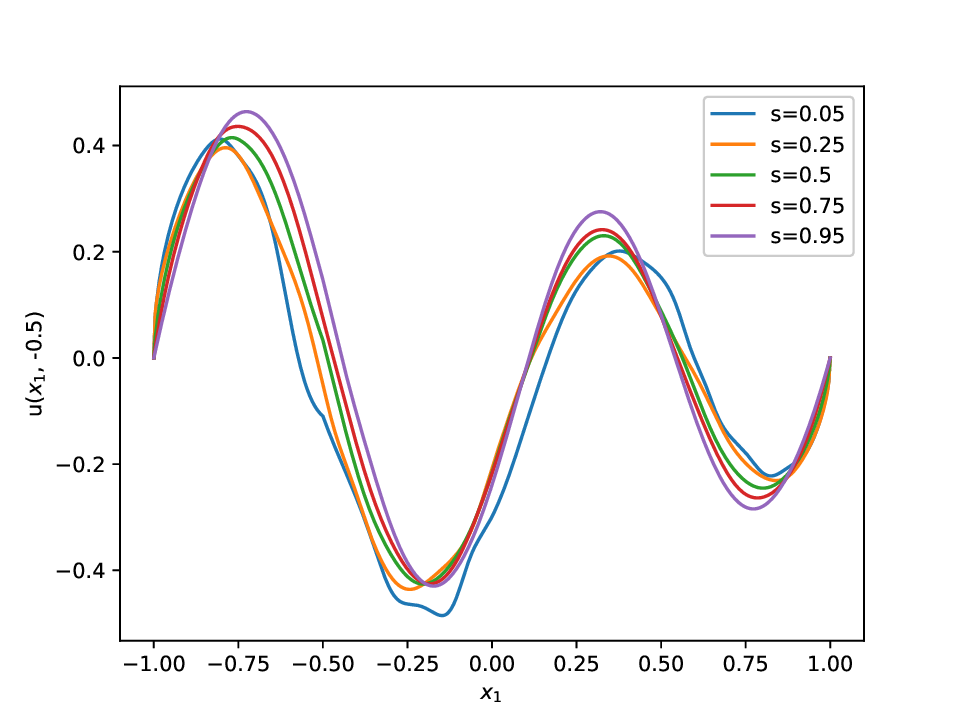}

{\small 
\parbox[t]{1.7cm}{\bf Figure 5}
\parbox[t]{10.5cm}{The eigenfunctions of \eqref{problem_laplace} in the $L$-shaped domain $[-1,1]^2 \backslash [0,1]^2$. 
	Each eigenfunction is normalized to let the maximum absolute value equal 1.
	The first column shows the $1$st eigenfunctions for different fractional orders $s$ while the second, third and fourth columns show the $5$th, $6$th and $10$th eigenfunctions respectively. 
	The first five rows show the eigenfunctions corresponding to $s = 0.05, 0.25, 0.5, 0.75$ and $0.95$ respectively. 
	The last row shows eigenfunctions at $x_2 = -0.5$.
}}

\end{center}
\end{figure}

\section{Isospectral problem}
\label{Section_Isospectral}

In this section, we explore the fractional order isospectral problem. In 1966, Kac posed the famous isospectral problem\supercite{drum1_kac1966can}, ``Can one hear the shape of a drum", which asks whether the Laplace operator with Dirichlet boundary conditions on two different domains can have the same spectrum?
In 1992, Gordan, Webb, and Wolpert gave a negative answer to this question with a counterexample\supercite{drum2_gordon1992one,drum2_gordon1992isospectral}, proving that it is possible for two domains to have the same spectrum.
Figure 6 represents their counterexample. 
Since then, researchers have discovered numerous pairs of domains with identical spectra, and the eigenvalues of many of them are calculated by some numerical works\supercite{drum3_driscoll1997eigenmodes, drum4_borzi2006algebraic, drum5_li2017efficient}.

Now, we wonder whether two different domains that have the same spectrum for the Laplace operator will also have the same spectrum for the fractional Laplace operator. We solve the eigenvalue problem of the fractional Laplace operator in these two domains to draw a conjecture to this question. 
The relative differences between two eigenvalues are calculated by
\begin{equation}
R_{k}^{(s)} = \frac{\lambda_{B,k}^{(s)} - \lambda_{A,k}^{(s)}}{\big(\lambda_{A,k}^{(s)} + \lambda_{B,k}^{(s)}\big)/2}.
\end{equation} 
Here, $\lambda_{A,k}^{(s)}$ and $\lambda_{B,k}^{(s)}$ are the $k$th eigenvalues for the fractional Laplacian with order $s$ in the domains of drum A and drum B, respectively.
The previous conclusion stated that 
\begin{equation}
	R_{k}^{(1)} = 0, \quad \text{for any k.}
\end{equation}
But based on the following experiments, we speculate that when $0<s<1$, 
\begin{equation}
	R_{k}^{(s)} \neq 0, \quad \text{for some $k$} .
\end{equation}

\begin{figure}
	\begin{center}
		\setcounter{subfigure}{0}
		\subfigure[Drum A]{
			\includegraphics[width=0.35\columnwidth]{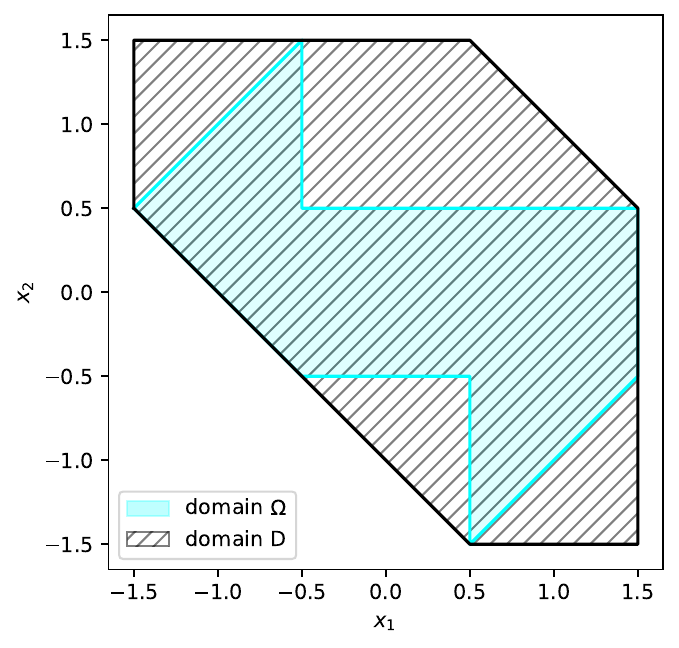}}
		\subfigure[Drum B]{
			\includegraphics[width=0.35\columnwidth]{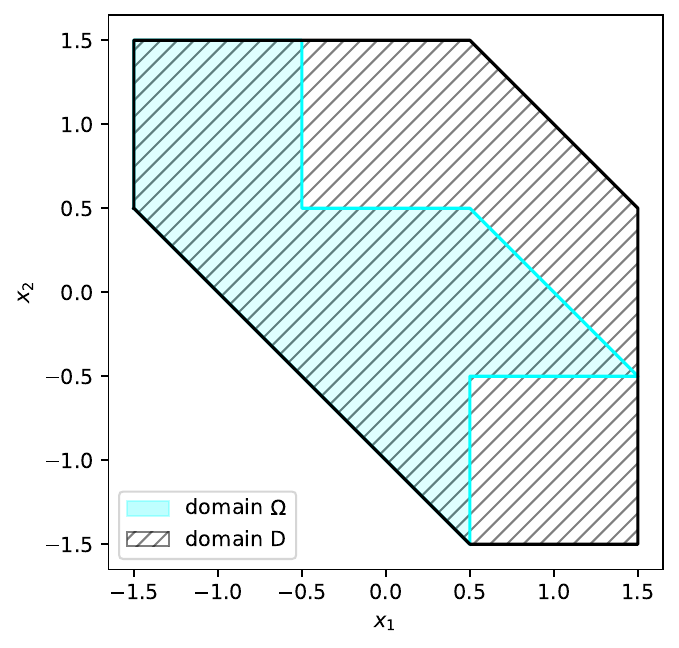}}
		
		\small {\bf Figure 6}\ \ Shape of the problem domain $\itOmega$ and the sampling domain $D$ of the drum-shaped problem.
		
	\end{center}
\end{figure}

Since these two domains are not convex, we select a convex domain $D$ for sampling and the shape of the domain $D$ is plotted in Figure 6 also. We construct two types of feature functions similar to the previous example. The first-type feature functions are used to capture the singularity near the boundary while the second-type captures the singularity at the corners. The network has a width of $60$ and consists of $40$ first-type feature functions and $20$ second-type feature functions. All other settings remain the same as the previous examples. 

Table 10 reports the first two eigenvalues in these two domains. It is evident that the first eigenvalue in the domain of drum A is smaller than that in drum B, whereas the second eigenvalue in drum A is greater than that in drum B.

We further compute more eigenvalues and calculate the relative difference between them. The results for different $k$ and $s$ are shown in Figure 7. The value $R_{k}^{(s)}$ for these $k$ and $s$ we calculated are significantly different from $0$ and the maximum relative difference reaches $1.8\%$. These discrepancies cannot be explained solely by the sampling error of the Monte Carlo method. Therefore, we conjecture that even if the spectra of two domains are identical when $s=1$, they would not be the same for $0<s<1$.

\begin{table}
\begin{center}
{\small

{\bf Table 10}\ \  Estimates of the first two eigenvalues of \eqref{problem_laplace} in two drum-shaped domains.
\vskip 1mm

\begin{tabular}{ccccccccccc}
	\hline
	&	s    & 0.1    & 0.2    & 0.3    & 0.4    & 0.5    & 0.6    & 0.7    & 0.8    & 0.9      \\ \hline
	\multirow{2}{*}{k=1} & Drum A & 1.1429 & 1.3406 & 1.6114 & 1.9815 & 2.4887 & 3.1880 & 4.1578 & 5.5182 & 7.4383 \\
	& Drum B & 1.1438 & 1.3437 & 1.6172 & 1.9921 & 2.5054 & 3.2131 & 4.1913 & 5.5537 & 7.4694\\ \hline
	\multirow{2}{*}{k=2} & Drum A & 1.2258 & 1.5244 & 1.9237 & 2.4663 & 3.2077 & 4.2342 & 5.6653 & 7.6851 & 10.566  \\
	& Drum B & 1.2222 & 1.5152 & 1.9056 & 2.4335 & 3.1559 & 4.1591 & 5.5648 & 7.5679 & 10.450 \\ \hline
\end{tabular}

}
\end{center}
\label{Table10}
\end{table}

\begin{figure}
	\begin{center}
		
		\includegraphics[width=0.75\columnwidth]{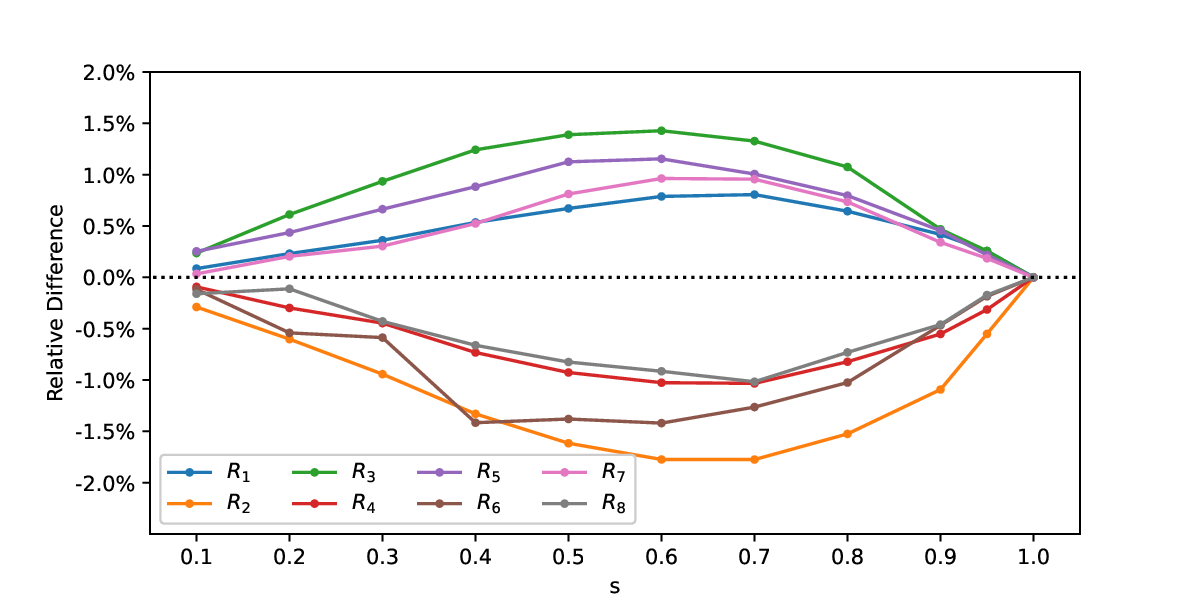}
		
		\small {\bf Figure 7}\ \ The relative difference $R_{k}^{(s)}$ for difference $k$ and $s$.
		
	\end{center}
\end{figure}
\section{Conclusion}\label{Section_Conclusion}
Based on a new loss function and a knowledge-based neural network architecture, we propose a novel deep learning method for computing eigenvalues of the fractional Schrödinger operator. 
We apply the method to problems in high-dimensional space and irregular domains in low dimensions. 
The numerical results demonstrate that the accuracy of our method in calculating the first few dozen eigenvalues of various problems, and this method outperforms the finite element method\supercite{FEM_borthagaray2018finite}.
We also draw a new conjecture to the fractional order isospectral problem for exhibiting the capability of the method.



\end{document}

%% file: jsscN.tex
\makeatletter
\def\@oddfoot{\hfill}
\newcount\shumeicount
\def\setshumei#1#2#3{%
  \shumeicount=\count0
  \def\@oddhead{%
    \raise-5pt\hbox to0pt{\vrule width\hsize height 0pt depth 0.4pt\hss}\relax
    \ifnum \shumeicount=\count0
      \raise-7pt\hbox to0pt{\vrule width\hsize height 0pt depth 0.4pt\hss}\relax
      #1
    \else
      \ifodd\count0
        #2
      \else
        #3
       \fi
     \fi
  }%
}
\makeatother
\makeatletter
\def\@oddfoot{\hfill}
\newcount\shujiaocount
\def\setshujiao{%
  \shujiaocount=\count0
  \def\@oddfoot{%
      \ifodd\count0
      \else
      \fi
  }%
}
\makeatother
\def\title#1#2#3#4{{
  \vspace*{0.3cm}
  \begin{flushleft} \Large\bf #1\end{flushleft}
  \vspace*{-0.2cm}
      \begin{flushleft}
      \bf #2
      \end{flushleft}
      \footnotetext{\hspace{-6mm} #3\\ #4}}}

\def\dshm#1#2#3#4
{\setshumei{
		{ }\hfill}
	{\hfill {\small #3}\hfill\hbox to0pt{\hss\thepage}}
	{\hbox to0pt{\thepage\hss}\hfill {\small #4}\hfill
	}
	\setshujiao}
\def\drd#1#2
{{\vskip 1cm\small \begin{flushleft}
 #1 \\
 #2\\
\end{flushleft}}}

\def\hat{\widehat}
\def\tilde{\widetilde}
\def\bar{\overline}
\def\epsilon{\varepsilon}


%% file: eigenvalue-arxiv.bbl
\begin{thebibliography}{99}

\bibitem{intro_ad1_metzler2000random} Metzler R and Klafter J, 
The random walk's guide to anomalous diffusion: a fractional dynamics approach, 
{\it Physics Reports}, 2000, {\bf 339}(1): 1--77.

\bibitem{intro_tf1_shlesinger1987levy} Shlesinger M F, West B J, and Klafter J,
L{\'e}vy dynamics of enhanced diffusion: Application to turbulence,
{\it Physical Review Letters}, 1987, {\bf 58}(11): 1100.

\bibitem{intro_pmf1_de2011fractional} de Pablo A, Quir{\'o}s F, Rodr{\'i}guez A, and V{\'a}zquez J L,
A fractional porous medium equation,
{\it Advances in Mathematics}, 2011, {\bf 226}(2): 1378--1409.

\bibitem{intro_qm1_laskin2000fractional} Laskin N,
Fractional quantum mechanics and l{\'e}vy path integrals,
{\it Physics Letters A}, 2000, {\bf 268}(4-6): 298--305.

\bibitem{intro_qm2_laskin2002fractional} Laskin N,
Fractional Schrödinger equation,
{\it Physical Review E}, 2002, {\bf 66}(5): 056108.

\bibitem{intro_qm3_longhi2015fractional} Longhi S,
Fractional Schrödinger equation in optics,
{\it Optics Letters}, 2015, {\bf 40}(6): 1117--1120.

\bibitem{intro_qm4_zhang2016pt} Zhang Y, Zhong H, Belić M R, Zhu Y, Zhong W, Zhang Y, Christodoulides D N, and Xiao M,
PT symmetry in a fractional Schrödinger equation,
{\it Laser \& Photonics Reviews}, 2016, {\bf 10}(3): 526--531.

\bibitem{intro_AFEM1_ainsworth2017aspects} Ainsworth M and Glusa C,
Aspects of an adaptive finite element method for the fractional Laplacian: a priori and a posteriori error estimates, efficient implementation and multigrid solver,
{\it Computer Methods in Applied Mechanics and Engineering}, 2017, {\bf 327}: 4--35.

\bibitem{intro_AFEM2_ainsworth2018towards} Ainsworth M and Glusa C,
Towards an efficient finite element method for the integral fractional Laplacian on polygonal domains, 
{\it Contemporary Computational Mathematics - A Celebration of the 80th Birthday of Ian Sloan}, Dick J, Kuo F Y, Woźniakowski H, Springer, 2018: 17--57.

\bibitem{intro_FD2_del2018robust} Del Teso F, Endal J, and Jakobsen E R,
Robust numerical methods for nonlocal (and local) equations of porous medium type. Part II: Schemes and experiments,
{\it SIAM Journal on Numerical Analysis}, 2018, {\bf 56}(6): 3611--3647.

\bibitem{intro_Spectral1_mao2016efficient} Mao Z, Chen S, and Shen J,
Efficient and accurate spectral method using generalized Jacobi functions for solving Riesz fractional differential equations,
{\it Applied Numerical Mathematics}, 2016, {\bf 106}: 165--181.

\bibitem{intro_Spectral2_xu2018spectral} Xu K and Darve E,
Spectral method for the fractional Laplacian in 2D and 3D,
{\it arXiv preprint}, 2018, arXiv:1812.08325.

\bibitem{intro_WOS1_kyprianou2018unbiased} Kyprianou A E, Osojnik A, and Shardlow T,
Unbiased 'walk-on-spheres' Monte Carlo methods for the fractional Laplacian,
{\it IMA Journal of Numerical Analysis}, 2018, {\bf 38}(3): 1550--1578.

\bibitem{intro_WOS2_shardlow2019walk} Shardlow T,
A walk outside spheres for the fractional Laplacian: fields and first eigenvalue,
{\it Mathematics of Computation}, 2019, {\bf 88}(320): 2767--2792.

\bibitem{intro_WOS3_sheng2022efficient} Sheng C, Su B, and Xu C,
Efficient Monte Carlo method for integral fractional Laplacian in multiple dimensions,
{\it arXiv preprint}, 2022, arXiv:2204.08860.



\bibitem{du2020acta} D'Elia M, Du Q, Glusa C, Gunzburger M, Tian X, and Zhou Z,
Numerical methods for nonlocal and fractional models,
{\it Acta Numerica}, 2020, {\bf 29}: 1--124.	

\bibitem{Karniadakis2020whatis} Lischke A, Pang G, Gulian M, Song F, Glusa C, Zheng X, Mao Z, Cai W, Meerschaert M M, Ainsworth M and Karniadakis G E,
What is the fractional Laplacian? A comparative review with new results,
{\it Journal of Computational Physics}, 2020, {\bf 404}: 109009.

\bibitem{intro_num_method3_bonito2018numerical} Bonito A, Borthagaray J P, Nochetto R H, Ot{\'a}rola E, and Salgado A J,
Numerical methods for fractional diffusion,
{\it Computing and Visualization in Science}, 2018, {\bf 19}: 19--46.


\bibitem{dyda2017eigenvalues} Dyda B, Kuznetsov A, and Kwa{\'s}nicki M,
Eigenvalues of the fractional Laplace operator in the unit ball,
{\it Journal of the London Mathematical Society}, 2017, {\bf 95}(2): 500--518.

\bibitem{dyda2012fractional} Dyda B,
Fractional calculus for power functions and eigenvalues of the fractional Laplacian,
{\it Fractional Calculus and Applied Analysis}, 2012, {\bf 15}(4): 536--555.

\bibitem{bao2020jacobi} Bao W, Chen L, Jiang X, and Ma Y,
A Jacobi spectral method for computing eigenvalue gaps and their distribution statistics of the fractional Schrödinger operator, 
{\it Journal of Computational Physics}, 2020, {\bf 421}: 109733.

\bibitem{FEM_borthagaray2018finite} Borthagaray J P, Del Pezzo L M, and Mart{\'\i}nez S,
Finite element approximation for the fractional eigenvalue problem,
{\it Journal of Scientific Computing}, 2018, {\bf 77}(1): 308--329.

\bibitem{intro_NN_early1_samardzija1991neural} Samardzija N and Waterland R L,
A neural network for computing eigenvectors and eigenvalues,
{\it Biological Cybernetics}, 1991, {\bf 65}(4): 211--214.

\bibitem{intro_NN_early2_cichocki1992neural} Cichocki A and Unbehauen R,
Neural networks for computing eigenvalues and eigenvectors,
{\it Biological Cybernetics}, 1992, {\bf 68}: 155--164.

\bibitem{intro_NN_mb1_carleo2017solving} Carleo G and Troyer M,
Solving the quantum many-body problem with artificial neural networks,
{\it Science}, 2017, {\bf 355}(6325), 602--606.

\bibitem{intro_NN_mb2_choo2020fermionic} Choo K, Mezzacapo A, and Carleo G, 
Fermionic neural-network states for ab-initio electronic structure,
{\it Nature Communications}, 2020, {\bf 11}(1): 2368.

\bibitem{intro_NN_mb5_han2019solving} Han J, Zhang L, and E W,
Solving many-electron Schrödinger equation using deep neural networks,
{\it Journal of Computational Physics}, 2019, {\bf 399}: 108929.

\bibitem{intro_NN_eigen1_han2020solving} Han J, Lu J, and Zhou M,
Solving high-dimensional eigenvalue problems using deep neural networks: A diffusion Monte Carlo like approach,
{\it Journal of Computational Physics}, 2020, {\bf 423}: 109792.

\bibitem{intro_NN_eigen2_li2022semigroup} Li H and Ying L,
A semigroup method for high dimensional elliptic PDEs and eigenvalue problems based on neural networks,
{\it Journal of Computational Physics}, 2022, {\bf 453}: 110939.

\bibitem{intro_NN_eigen3_simonnet2022deep} Simonnet E and Chekroun M D,
Deep spectral computations in linear and nonlinear diffusion problems,
{\it arXiv preprint}, 2022, arXiv:2207.03166.

\bibitem{intro_NN_eigen5_zhang2022solving} Zhang W, Li T, and Schütte C,
Solving eigenvalue PDEs of metastable diffusion processes using artificial neural networks, {\it Journal of Computational Physics}, 2022, {\bf 465}: 111377.

\bibitem{intro_NN_eigen9_elhamod2022cophy} Elhamod M, Bu J, Singh C, Redell M, Ghosh A, Podolskiy V, Lee W, and Karpatne A, CoPhy-PGNN: Learning physics-guided neural networks with competing loss functions for solving eigenvalue problems, 
{\it ACM Transactions on Intelligent Systems and Technology}, 2022, {\bf 13}(6): 1--23.

\bibitem{intro_NN_eigen10_finol2019deep} Finol D, Lu Y, Mahadevan V, and Ankit S,
Deep convolutional neural networks for eigenvalue problems in mechanics,
{\it International Journal for Numerical Methods in Engineering}, 2019, {\bf 118}(5): 258--275.

\bibitem{intro_NN_eigen_theory_lu2022priori} Lu J and Lu Y,
A priori generalization error analysis of two-layer neural networks for solving high dimensional Schrödinger eigenvalue problems, 
{\it Communications of the American Mathematical Society}, 2022, {\bf 2}(01): 1--21.

\bibitem{kwasnicki2017ten} Kwa{\'s}nicki M,
Ten equivalent definitions of the fractional Laplace operator, 
{\it Fractional Calculus and Applied Analysis}, 2017, {\bf 20}(1): 7--51.

\bibitem{valdinoci2009long} Valdinoci E,
From the long jump random walk to the fractional Laplacian,
{\it Boletín de la Sociedad Española de Matemática Aplicada}, 2009, {\bf 49}: 33--44.



\bibitem{property1_frank2016eigenvalue} Frank R L, 
Eigenvalue Bounds for the Fractional Laplacian: A Review, 
{\it Recent Developments in Nonlocal Theory}, Palatucci G and Kuusi T, De Gruyter, 2017: 210--235.

\bibitem{property2_chen2005two} Chen Z Q and Song R,
Two-sided eigenvalue estimates for subordinate processes in domains,
{\it Journal of Functional Analysis}, 2005, {\bf 226}(1): 90--113.

\bibitem{property3_bao2018fundamental} Bao W, Ruan X, Shen J, and Sheng C,
Fundamental gaps of the fractional Schr{\"o}dinger operator,
{\it Communications in Mathematical Sciences}, 2019, {\bf 17}(2): 447--471.


\bibitem{boundary_order_grubb2015fractional} Grubb G,
Fractional Laplacians on domains, a development of H{\"o}rmander's theory of $\mu$-transmission pseudodifferential operators,
{\it Advances in Mathematics}, 2015, {\bf 268}: 478--528.

\bibitem{architecture_khoo2019solving} Khoo Y, Lu J, and Ying L.
Solving for high-dimensional committor functions using artificial neural networks,
{\it Research in the Mathematical Sciences}, 2019, {\bf 6}: 1-13.

\bibitem{Laplacian_eigenvalue_square3_liu2013verified} Liu X and Oishi S,
Verified eigenvalue evaluation for the Laplacian over polygonal domains of arbitrary shape, 
{\it SIAM Journal on Numerical Analysis}, 2013, {\bf 51}(3): 1634--1654.

\bibitem{Laplacian_eigenvalue_Lshaped2_yuan2009bounds} Yuan Q and He Z,
Bounds to eigenvalues of the Laplacian on L-shaped domain by variational methods,
{\it Journal of Computational and Applied Mathematics}, 2009, {\bf 233}(4): 1083--1090.

\bibitem{drum1_kac1966can} Kac M,
Can one hear the shape of a drum?,
{\it The American Mathematical Monthly}, 1966, {\bf 73}(492):1--33.

\bibitem{drum2_gordon1992one} Gordon C, Webb D L, and Wolpert S,
One cannot hear the shape of a drum,
{\it Bulletin of the American Mathematical Society}, 1992, {\bf 27}(1): 134--138.

\bibitem{drum2_gordon1992isospectral} Gordon C, Webb D L, and Wolpert S,
Isospectral plane domains and surfaces via Riemannian orbifolds,
{\it Inventiones Mathematicae}, 1992, {\bf 110}(1): 1--22.

\bibitem{drum3_driscoll1997eigenmodes} Driscoll T A,
Eigenmodes of isospectral drums,
{\it SIAM Review}, 1997, {\bf 39}(1): 1--17.

\bibitem{drum4_borzi2006algebraic} Borzì A and Borzì G,
Algebraic multigrid methods for solving generalized eigenvalue problems,
{\it International Journal for Numerical Methods in Engineering}, 2006, {\bf 65}(8): 1186--1196.

\bibitem{drum5_li2017efficient} Li H and Zhang Z,
Efficient spectral and spectral element methods for eigenvalue problems of Schrödinger equations with an inverse square potential,
{\it SIAM Journal on Scientific Computing}, 2017, {\bf 39}(1): A114--A140.
  
\end{thebibliography}
